\documentclass[11pt]{amsart}
\usepackage{amsmath, amssymb, amsthm, amsfonts}
\usepackage[all]{xy}
\usepackage{graphicx}
\usepackage{subcaption}
\usepackage{tabularx}
\usepackage{hyperref}
\usepackage[utf8]{inputenc}

\usepackage{float}
\usepackage{enumerate}

\theoremstyle{plain}
\newtheorem{Thm}{Theorem}[section]
\newtheorem{Prop}[Thm]{Proposition}
\newtheorem{Cor}[Thm]{Corollary}
\newtheorem{Lem}[Thm]{Lemma}
\newtheorem{Conj}[Thm]{Conjecture}

\theoremstyle{definition}
\newtheorem{Def}[Thm]{Definition}

\theoremstyle{remark}
\newtheorem{Rmk}[Thm]{Remark}
\newtheorem{Ex}[Thm]{Example}

\begin{document}

\title[A matryoshka structure of higher secant varieties]{A matryoshka structure of higher secant varieties and the generalized Bronowski's conjecture}
\author[Junho Choe and Sijong Kwak]{Junho Choe and Sijong Kwak
}
\address{Junho Choe \\
Department of Mathematics, Korea Advanced Institute of Science and Technology (KAIST), 373-1 Gusung-dong, Yusung-Gu, Daejeon, Republic of Korea}
\email{junhochoe@kaist.ac.kr}
\address{Sijong Kwak\\
Department of Mathematics, Korea Advanced Institute of Science and Technology (KAIST), 373-1 Gusung-dong, Yusung-Gu, Daejeon, Republic of Korea}
\email{sjkwak@kaist.ac.kr}


\begin{abstract}
In projective algebraic geometry, there are classical and fundamental results that describe the structure of geometry and syzygies, and many of them characterize varieties of minimal degree and del Pezzo varieties. In this paper, we consider analogous objects in the category of higher secant varieties. Our main theorems say that there is a matryoshka structure among those basic objects including a generalized $K_{p,1}$ theorem, syzygetic and geometric characterizations of higher secant varieties of minimal degree and del Pezzo higher secant varieties, defined in this paper.
For our purpose, we prove a weak form of the generalized Bronowski's conjecture raised by C. Ciliberto and F. Russo that relates the identifiability for higher secant varieties to the geometry of tangential projections.
\end{abstract}

\maketitle
\tableofcontents 

\section{Introduction}

In classical projective algebraic geometry, varieties of minimal degree and del Pezzo varieties are two basic objects, and have been classified and characterized in many aspects. Such viewpoints provide structural results on geometry and syzygies. Among them, we consider the following themes:
\begin{itemize}
\item Castelnuovo's and Fano's results on the number of quadrics (\cite{castelnuovo1889ricerche} and \cite{fano1893sopra});

\item M. Green's $K_{p,1}$ theorem (\cite{green1984koszul});

\item Eisenbud-Goto characterization of $2$-regular projective varieties (\cite{eisenbud1984linear});

\item The rigidity result on property $N_p$ (\cite{eisenbud2005restricting}); and


\item Projective varieties of almost minimal degree.
\end{itemize}
Some of the subjects above are classical work, and the others are somewhat recent work. These have been studied by many authors, and developed in various directions for a long time.
 
  On the other hand, {\it higher secant varieties} are geometric objects studied for over a hundred years in algebraic geometry. Historically, much attention has been paid to the defectiveness of higher secant varieties due to F. Severi, F. Zak, L. Chiantini and C. Ciliberto, etc
(\cite{chiantini2010dimension}, \cite{severi1901intorno} and \cite{zak1983projections}).
Recently, computing the exact dimension of higher secant varieties to specific varieties has also been important in tensor geometry and algebraic variants of P versus NP problem in complexity theory (\cite{landsberg2012tensors}).

In this paper, we are interested in the category of higher secant varieties, and revisit the aforementioned interplay between geometry and syzygies, especially from a viewpoint of {\it matryoshka structure}. 

``Matryoshka" is a set of traditional Russian dolls with the same shape but different sizes inside successively. Mathematically, many fundamental results on geometry and syzygies for projective varieties are expected to hold in the category of higher secant varieties in a very similar pattern. In 2015, F. Zak pointed out to the second author this kind of phenomena implicitly to which he gave the name ``matryoshka structure.'' 

The authors have realized that structures of the Betti tables and their geometric interpretations in the category of higher secant varieties are very similar to those in the category of projective varieties just like ``Matryoshka". In fact, we are able to find counterparts of previous results listed above for higher secant varieties. Surprisingly, it was naturally done to characterize interesting higher secant varieties acting like del Pezzo varieties. (Another kind of matryoshka phenomenon was also mentioned in the book \cite{russo2016geometry}. It concerns the comparison of uniruled varieties and their varieties of lines, and provides surprising applications to the classification of Severi varieties and Hartshorne's conjecture on complete intersections.)

We always work over an algebraically closed field of characteristic zero, and let $X$ be a nondegenerate irreducible reduced projective variety in $\mathbb P^r$. For any positive integer $q\ge 1$, the {\it $q$-secant variety}, denoted by $S^q(X)$, is the Zariski closure of the union of all linear spaces of dimension $q-1$ spanned by linearly independent $q$ points of $X$. It was proved in \cite{ciliberto2006varieties} that
$$
\deg S^q(X)\geq\binom{e+q}{q},
$$
where $e$ is the codimension of $S^q(X)$ in $\mathbb P^r$. One can say that $S^q(X)$ is a {\it $q$-secant variety of minimal degree} if $\deg S^q(X)=\binom{e+q}{q}$, and $X$ is called an {\it $\textup{M}^q$-variety} if $S^q(X)$ is a $q$-secant variety of minimal degree. Thus, for each $q\ge 1$, there are objects of minimal degree in the category of $q$-secant varieties, and it would be very interesting to characterize or classify objects of minimal degree: For $q=1$, there are many characterization results on minimal degree varieties in various aspects, and it is well known that varieties of minimal degree are only quadric hypersurfaces, a cone over the Veronese surface in $\mathbb P^5$ and rational normal scrolls (\cite{aprodu2010koszul}, \cite{eisenbud1987varieties}, \cite{green1984koszul}). On the other hand, for $q\ge 2$, it would be very difficult to do the same thing in the category of $q$-secant varieties . Note that for a variety $X$ of minimal degree, $S^q(X)$ is also of minimal degree (\cite[Claim 5.2]{ciliberto2006varieties} or Corollary \ref{minmin}). But the converse is not true in general. For irreducible curves, the converse is also true. For smooth surfaces, some exceptions were classified due to Ciliberto-Russo (\cite[Theorem 9.2]{ciliberto2006varieties}), and however, it might be a very hard task to keep going on this kind of classification.

 Instead, through detailed observation of the matryoshka structure, we establish the following syzygetic characterization:

\begin{Thm}\label{syzM}
Let $e=\textup{codim}S^q(X)$. Assume that $e\geq 1$ and $q\ge 1$. Then, we have the following equivalence:
\begin{enumerate}[\normalfont (1)]
\item[\normalfont (1)] $S^q(X)\subset\mathbb P^r$ is a $q$-secant variety of minimal degree.
\item[\normalfont (2)] $S^q(X)\subset\mathbb P^r$ has $q$-pure Cohen-Macaulay Betti table.
\item[\normalfont (3a)] $\beta_{p,q}(S^q(X))=\binom{p+q-1}{q}\binom{e+q}{p+q}$ for every $p$.
\item[\normalfont (3b)] $\dim(I_{S^q(X)})_{q+1}=\binom{e+q}{q+1}$.
\item[\normalfont (3c)] $\beta_{p,q}(S^q(X))=\binom{p+q-1}{q}\binom{e+q}{p+q}$ for some $1\leq p\leq e$.
\item[\normalfont (4)] the $(q+1)$-strand of $S^q(X)\subset\mathbb P^r$ has length $e$.
\item[\normalfont (5a)] $\operatorname{reg}S^q(X)=q+1$.
\item[\normalfont (5b)] the homogeneous coordinate ring $S_{S^q(X)}$ satisfies property $N_{q+1,e}$.
\end{enumerate}
\end{Thm}
\noindent For instance, (2) tells us that the Betti table of any minimal degree $q$-secant variety is the following:
\begin{figure}[H]
\centering
\begin{tabularx}{0.5\linewidth}{ c| *{5}{>{\centering\arraybackslash}X} }
       & $0$ & $1$ & $\cdots$ & $e-1$ & $e$ \\ \hline
$0$ & $1$ & & & & \\
$q$ & & $\ast$ & $\cdots$ & $\ast$ & $\ast$
\end{tabularx}
\caption{The Betti table of $q$-secant varieties of minimal degree}
\end{figure}
\noindent Refer to the fact that there are no hypersurfaces of degree $q$ containing $S^q(X)$.

In fact, it is proved in Corollary \ref{syzbound} that $\beta_{p,q}(S^q(X))\leq\binom{p+q-1}{q}\binom{e+q}{p+q}$ for every $p$. Thus, among all $q$-secant varieties, $q$-secant varieties of minimal degree are characterized  to have the maximal Betti number $\beta_{p,q}(S^q(X))=\binom{p+q-1}{q}\binom{e+q}{p+q}$ for some $1\leq p\leq e$. For $q=1$, it is classical and was proved by G. Castelnuovo, Eisenbud-Green-Hulek-Popescu and Han-Kwak (\cite {castelnuovo1889ricerche},\cite{eisenbud2005restricting},\cite {han2015sharp}).

On the other hand, in birational geometry, Fano varieties including del Pezzo surfaces are basic objects as building blocks. However, together with varieties of minimal degree, del Pezzo varieties are basic objects in projective algebraic geometry, where by definition, a del Pezzo variety is an arithmetically Cohen-Macaulay variety of almost minimal degree. One may expect that there is an analogue of the del Pezzo variety in the category of higher secant varieties. Note that varieties of almost minimal degree have sectional genus less than or equal to one, and that a variety of almost minimal degree has sectional genus one if and only if it is a del Pezzo variety. Similarly, suppose that $S^q(X)$ is not a $q$-secant variety of minimal degree with codimension $e\geq 1$. Then, it can be shown (Proposition \ref{degboundII}) that
$$
\deg S^q(X)\geq\binom{e+q}{q}+\binom{e+q-1}{q-1}.
$$
In addition, if $\deg S^q(X)=\binom{e+q}{q}+\binom{e+q-1}{q-1}$, then we obtain
the upper bound of the sectional genus (Proposition \ref{secgenbound}):
$$
\pi(S^q(X))\leq(q-1)\left(\binom{e+q}{q}+\binom{e+q-1}{q-1}\right)+1.
$$
(Note that in positive characteristic, the degree bound above is not true in general (see Remark \ref{positive}).) It would be natural to call $S^q(X)$ a {\it del Pezzo $q$-secant variety} if $\deg S^q(X)=\binom{e+q}{q}+\binom{e+q-1}{q-1}$ and $\pi(S^q(X))=(q-1)\left(\binom{e+q}{q}+\binom{e+q-1}{q-1}\right)+1$. This definition coincides with the equivalent definition of del Pezzo varieties. The existence of del Pezzo higher secant varieties can be expected in view of the matryoshka structure.

As in the case of minimal degree $q$-secant varieties, the viewpoint of the matryoshka structure enables us to prove the following:
\begin{Thm}\label{syzdP}
Put $e=\textup{codim}S^q(X)$. Assume that $e\geq 2$ and $q\ge 1$. The following are equivalent:
\begin{enumerate}[\normalfont (1)]
\item[\normalfont (1)] $S^q(X)\subset\mathbb P^r$ is a del Pezzo $q$-secant variety.

\item[\normalfont (2)] $S^q(X)\subset\mathbb P^r$ has $q$-pure Gorenstein Betti table.

\item[\normalfont (3a)] $\beta_{p,q}(S^q(X))=\binom{p+q-1}{q}\binom{e+q}{p+q}-\binom{e+q-p-1}{q-1}\binom{e+q-1}{e+q-p}$ for every $1\leq p\leq e$.

\item[\normalfont (3b)] $\dim(I_{S^q(X)})_{q+1}=\binom{e+q}{q+1}-\binom{e+q-2}{q-1}$.

\item[\normalfont (3c)] $\beta_{p,q}(S^q(X))=\binom{p+q-1}{q}\binom{e+q}{p+q}-\binom{e+q-p-1}{q-1}\binom{e+q-1}{e+q-p}$ for some $1\leq p\leq e-1$.

\item[\normalfont (4)] $S_{S^q(X)}$ satisfies property $N_{q+1,e-1}$, but not property $N_{q+1,e}$.
\end{enumerate}
\end{Thm}
\noindent Hence, every del Pezzo $q$-secant variety has Betti table of the following form:
\begin{figure}[H]
\centering
$$
\begin{tabularx}{0.5\linewidth}{ c| *{5}{>{\centering\arraybackslash}X} }
         & $0$ & $1$     & $\cdots$ & $e-1$ & $e$ \\ \hline
$0$   & $1$ &            &               &          & \\
$q$   &        & $\ast$ & $\cdots$ & $\ast$ & \\
$2q$ &        &            &             &           & $1$
\end{tabularx}
$$
\caption{The Betti table of del Pezzo $q$-secant varieties}
\end{figure}

For $q=1$, it is classical and was proved by G. Fano, and Han-Kwak (\cite{fano1893sopra}, \cite {han2015sharp}). It was also proved for higher secant varieties to elliptic normal curves due to Bothmer-Hulek, or T. Fisher (see \cite{bothmer2004geometric},\cite{fisher2006higher}).
It looks very mysterious that $S^q(X)\subset\mathbb P^r$ is not del Pezzo $q$-secant even if $X$ is del Pezzo; for smooth del Pezzo surfaces, $S^q(X)$ is always a minimal degree $q$-secant variety. However, for curves, we prove that $X$ is del Pezzo if and only if $S^q(X)$ is also a del Pezzo $q$-secant variety (Proposition \ref{gBcIIcurve}).

In \cite{bronowski1933sum}, J. Bronowski studied the Waring problem for homogeneous polynomials, and he claimed that when $n+1$ divides $\binom{d+n}{n}$ with quotient $q$, a general homogeneous form of degree $d$ in $n+1$ variables can be expressed as the sum of $q$ linear forms to the power of $d$ if and only if a general $(q-1)$-tangential projection of the Veronese variety $\nu_d(\mathbb P^n)\subseteq\mathbb P^{\binom{d+n}{n}-1}$ is birational, and that a similar thing holds for arbitrary (smooth) projective varieties. However, the proof was considered insufficient by Ciliberto-Russo. They extended the Bronowski's assertion, and called it the {\it generalized Bronowski's conjecture}.

One says $X$ to be a {\it variety with the minimal number of apparent $q$-secant $(q-2)$-planes}, or an {\it $\textup{MA}^q$-variety} if it is an $\text{M}^q$-variety and a unique span of linearly independent $q$ points of $X$ passes through a general point of $S^q(X)$, namely $X$ is {\it identifiable} for $S^q(X)$. As noted in \cite{ciliberto2006varieties}, $X^n\subset\mathbb P^{2n+1}$ is an $\text{MA}^2$-variety if and only if it is a {\it variety with one apparent double point} which was studied by Severi. For an $\text{MA}^q$-variety $X$, its general $(q-1)$-tangential projection is birational onto a variety of minimal degree. Its proof depends on the inequality that the number of $q$-secant $(q-1)$-planes to $X$ passing through a general point of $S^q(X)$ is at least $\deg\pi_\Lambda\deg X_\Lambda/\deg S^q(X)$ when $\dim S^q(X)=q\dim X+q-1$ with $\pi_\Lambda$ and $X_\Lambda$ a general $(q-1)$-tangential projection of $X$; see \cite[Theorem 2.7 and Corollary 4.5]{ciliberto2006varieties}.
\begin{Conj}[The generalized Bronowski's conjecture, Ciliberto-Russo]
If a general $(q-1)$-tangential projection of $X$ is birational and its image is a variety of minimal degree, then $X$ is an $\textup{MA}^q$-variety.
 \end{Conj}

To obtain our main theorems on characterization, we need to prove a weak form of the generalized Bronowski's conjecture for $q$-secant varieties of minimal degree which was described for the first time due to Ciliberto-Russo, and also to prove an analogue of the weak generalized Bronowski's conjecture for del Pezzo higher secant varieties.

\begin{Thm}[The generalized Bronowski's conjecture, weak form]\label{gBcwf}
Let $e$ be the codimension of $S^q(X)$. If either
\begin{enumerate}[\normalfont (1)]
\item a general inner projection $X_z$ is an $\textup{M}^q$-variety for $e\geq 2,q\geq 1$; or
\item a general tangential projection $X_{\mathbb T_zX}$ is an $\textup{M}^{q-1}$-variety for $e\geq 1,q\geq 2$,
\end{enumerate}
then $X$ is an $\textup{M}^q$-variety.
\end{Thm} 
\noindent It is elementary to check that any $\text{M}^q$-variety satisfies both (1) and (2). Consequently, when $S^q(X)\neq\mathbb P^r$, $X$ is an $\text{M}^q$-variety if and only if a general $(q-1)$-tangential projection of $X$ is a variety of minimal degree (Corollary \ref{gBcwfCor}). 

Once one shows that there passes a unique $(q-1)$-plane spanned by $q$ points in $X$ through a general point of $S^q(X)$ under the assumption, the original conjecture is proved. This is the case for most of projective varieties including curves; we verify the generalized Bronowski's conjecture when either $X$ is a curve or satisfies $\dim S^{q+1}(X)=\dim\widetilde{S}^{q+1}(X)$:
\begin{Thm}[Corollary \ref{gBccurve} and Theorem \ref{gBcusual}]
For a projective variety $X$ and its $q$-secant variety $S^q(X)$, the generalized Bronowski's conjecture is true when either
\begin{enumerate}[\normalfont (1)]
\item $\dim S^{q+1}(X)=(q+1)\dim X+q$; or

\item $X$ is a curve.
\end{enumerate}
\end{Thm}
\noindent Thus, the Ciliberto-Russo conjecture is still open and interesting for the following cases: $S^{q+1}(X)$ is nondefective with $S^{q+1}(X)=\mathbb P^r$, or defective. 

And as an application of the weak form of the generalized Bronowski's conjecture, we extend M. Green's famous theorem, the {\it $K_{p,1}$ theorem}, to higher secant varieties. While the $K_{p,1}$ theorem gives an upper bound of the length of quadratic strand, and characterize its maximal and almost maximal cases, we establish an upper bound of the length of $(q+1)$-strand for $q$-secant varieties, and conclude that the corresponding maximal case is precisely $q$-secant varieties of minimal degree:
\begin{Thm}[A generalized $K_{p,1}$ theorem]\label{kp1}
Let $e:=\textup{codim}S^q(X)$. We have the following:
\begin{enumerate}[\normalfont (1)]
\item $K_{p,q}(S^q(X))=0$ whenever $p>e$.

\item $K_{e,q}(S^q(X))\neq 0$ if and only if $S^q(X)$ has minimal degree.
\end{enumerate}
\end{Thm}

On the other hand, let us call $X$ a {\it $\textup{dP}^q$-variety} if $S^q(X)$ is a del Pezzo $q$-secant variety, and define a $\text{dP}^q$-variety $X$ to be a {\it $\textup{dPA}^q$-variety} if it is identifiable for $S^q(X)$. As in the case of $\text{M}^q$-varieties, several parallel results hold: Being a $\text{dP}^q$-variety and a $\text{dPA}^q$-variety are both preserved under taking a general inner projection, a general tangential projection of a $\text{dP}^q$-variety becomes a $\text{dP}^{q-1}$-variety, and a general $(q-1)$-tangential projection of a $\text{dPA}^q$-variety is birational onto a del Pezzo variety (cf.\phantom{.}Remark \ref{dPfunda}(2) and Proposition \ref{gBcIIconverse}). 

As Ciliberto-Russo raised the generalized Bronowski's conjecture for $\text{MA}^q$-varieties, we may ask whether the converses of results mentioned above hold or not. In particular, we are interested in $\text{dPA}^q$-varieties and their general $(q-1)$-tangential projections. We would like to propose the {\it generalized Bronowski's conjecture II} (Conjecture \ref{gBcII}) which says that $X$ is a $\text{dPA}^q$-variety if a general $(q-1)$-tangential projection of $X$ is birational onto a del Pezzo variety of codimension two or more. We verify if for curves, and in addition, also have the following theorem for del Pezzo $q$-secant varieties:

\begin{Thm}[The generalized Bronowski's conjecture II, weak form]\label{gBcIIwf}
Let $e=\textup{codim}S^q(X)$. If either
\begin{enumerate}[\normalfont (1)]
\item a general inner projection $X_z$ is a $\textup{dP}^q$-variety for $e\geq 3,q\geq 1$; or
\item a general tangential projection $X_{\mathbb T_zX}$ is a $\textup{dP}^{q-1}$-variety for $e\geq 2,q\geq 3$,
\end{enumerate}
then $X$ is a $\textup{dP}^q$-variety.
\end{Thm}

 It would also be interesting to classify smooth projective surfaces that are $\text{dP}^q$-varieties for some integer $q$; they are the simplest case of $\text{dP}^q$-varieties except curves. Note that smooth del Pezzo surfaces do not give del Pezzo higher secant varieties.

\bigskip

The paper is organized as follows: In Section 2, we discuss inner and tangential projections of higher secant varieties and establish the basic structure of Betti tables for higher secant varieties. In Section 3, we prove a weak form of the generalized Bronowski's conjecture, and in Section 4, show the generalization of $K_{p,1}$ theorem for higher secant varieties. In Section 5, we investigate the algebraic structure 
of higher secant varieties of minimal degree. In Section 6, we define del Pezzo higher secant varieties with some basic background. Furthermore, the algebraic structures and syzygetic characterization of del Pezzo secant varieties are described. Finally, we end with remarks, open questions and an appendix for a systematic approach to the partial elimination ideals in order to study algebraic structures of higher secant varieties.

\bigskip

\noindent{\bf Acknowledgements.}
We would like to thank F. Zak for conversations on the matryoshka structure and interesting comments. We also thank D.\ Eisenbud and M.\ Mella for general remarks on higher secant varieties of minimal degree and the identifiability, and F.-O.\ Schreyer for helpful conversations on del Pezzo higher secant varieties. We are particularly grateful to F.\ Russo who read the first version of this paper and explained to us what he did on matryoshka phenomena. This work was supported by the National Research Foundation of Korea (NRF) grant funded by the Korea government (2017R1E1A1A03070765).

\section{Higher secant varieties}

\begin{Def}
Let $X\subseteq\mathbb P^r$ be an irreducible projective variety, and let $q\geq1$ be an integer. The {\it $q$-secant variety} $S^q(X)$ to $X$ is defined by
$$
S^q(X):=\overline{\bigcup_{z_i\in X}\langle z_1,\ldots,z_q\rangle}\subseteq\mathbb P^r\quad\text{with}\quad\langle z_1,\ldots,z_q\rangle=\mathbb P^{q-1},
$$
where the union runs over all linearly independent $q$ points $z_i$ in $X$.
\end{Def}

In addition, we define the {\it abstract $q$-secant variety} $\widetilde{S}^q(X)$ by
$$
\widetilde{S}^q(X)=\overline{\{(\zeta,w):\zeta\in U,w\in\langle\zeta\rangle\}}\subseteq X^{(q)}\times\mathbb P^r,
$$
where $X^{(q)}$ is the $q$-th symmetric power of $X$ , and $U$ is the set of linearly independent distinct $q$ points $\zeta\in X^{(q)}$. The second projection $X^{(q)}\times\mathbb P^r\to\mathbb P^r$ induces a surjective map
$$
\widetilde{S}^q(X)\to S^q(X),
$$
and this map will be frequently mentioned. If the map $\widetilde{S}^q(X)\to S^q(X)$ is generically finite of degree $d$, then a general point of $S^q(X)$ lies in exactly $d$ spans of linearly independent $q$ points of $X$, and vise versa.

In the category of $q$-secant varieties, we shall have many results that involve $(q-1)$-secant varieties. The parallel results for the case $q=1$ still hold subject to $S^0(X)=\emptyset\subset\mathbb P^r$, where the homogeneous ideal of the empty set $\emptyset$ is chosen to be the maximal homogeneous ideal so that $\text{codim}\emptyset=r+1$ and $\deg\emptyset=1$. For this reason, we define the $0$-secant variety $S^0(X)$ to be the empty set $\emptyset\subset\mathbb P^r$.

\bigskip

\noindent{\bf Notation.} From now on, we fix our notation:
\begin{itemize}
\item $X\subseteq\mathbb P^r$ is a nondegenerate irreducible reduced projective variety.

\item $\text{Vert}X$ is the vertex of $X$.

\item $S$ is the homogeneous coordinate ring of $\mathbb P^r$, and $V$ is the vector space $S_1$ of linear forms.

\item $I_X$ and $S_X$ are the homogeneous ideal and coordinate ring of $X\subseteq\mathbb P^r$, respectively.

\item $S_z$ is the polynomial ring in $S$ generated by linear forms vanishing at a point $z\in\mathbb P^r$, and $V_z$ is the vector space $(S_z)_1$.

\item $(I)$ is the ideal of $S$ generated by an ideal $I$ of a subring $R\subseteq S$.
\end{itemize}

\subsection{Projections and partial elimination ideals}

Given a $\lambda$-plane $\Lambda\subsetneq\mathbb P^r$, the {\it projection of $X$ from $\Lambda$} is the map $\pi_\Lambda:X\setminus\Lambda\to\mathbb P^{r-\lambda-1}$ defined by
$$
\pi_\Lambda(x_0:\cdots:x_\lambda:x_{\lambda+1}:\cdots:x_r)=(x_{\lambda+1}:\cdots:x_r)
$$
if $\Lambda=V(x_{\lambda+1},\ldots,x_r)$. By abuse of language, its image $X_\Lambda:=\overline{\pi_\Lambda(X\setminus\Lambda)}\subseteq\mathbb P^{r-\lambda-1}$ is also called a projection of $X$ from $\Lambda$.

By the irreducibility of $X$, the projection commutes with taking higher secant varieties, namely
$$
S^q(X)_\Lambda=S^q(X_\Lambda).
$$

\bigskip

\noindent{\bf Inner projections.}
When $\Lambda$ is a point of $X$, the projection is called an {\it inner projection}. The {\it $m$-inner projection} is to do inner projections $m$ times. We will frequently consider a general $m$-inner projection image $X_{\langle z_1,\ldots,z_m\rangle}$, where $z_1,\ldots,z_m$ are $m$ general points of $X$.

Let $z\in X$ be a point. Then, $\dim X-1\leq\dim X_z\leq\dim X$, and the right equality holds if and only if $z$ is not in $\text{Vert}X$. Choosing suitable coordinates $x_0,x_1,\ldots,x_r$, we may assume that $z=(1:0:\cdots:0)$. Then, the {\it projective tangent cone} $\mathbb TC_zX\subseteq\mathbb P^r$ to $X$ and the {\it projectivized tangent cone} $\mathbb PTC_zX\subseteq\mathbb P^{r-1}$ to $X$ are (possibly non-reduced) projective schemes defined by the leading coefficient $f'$ of $f=f'x_0^d+\text{(lower $x_0$-degree terms)}\in I_X$ in $x_0$ with $f'$ in variables $x_1,\ldots,x_r$. Note that $\mathbb TC_zX$ and $\mathbb PTC_zX$ are equidimensional with $\dim\mathbb TC_zX=\dim X$ and $\dim\mathbb PTC_zX=\dim X-1$. Whenever $z\not\in\text{Vert}X$, we have the degree formula
\begin{equation}\label{degformula}
\deg X=\deg \pi_z\deg X_z+\deg\mathbb TC_zX,
\end{equation}
where $\pi_z:X\dashrightarrow X_z$ is the projection map.

On the other hand, due to M. Green, the projection from one point can also be analyzed via the {\it partial elimination ideal} theory: We restrict to the inner projection case $z\in X$. Let $S_z$ be the homogeneous coordinate ring of $\mathbb P^{r-1}$, the image of projection map $\pi_z:\mathbb P^r\setminus\{z\}\to\mathbb P^{r-1}$. Then, there is an infinite filtration of $I_X$ by graded $S_z$-modules:
$$
\widetilde{K}_0(z,I_X)\subseteq\widetilde{K}_1(z,I_X)\subseteq\cdots\subseteq\widetilde{K}_i(z,I_X)\subseteq\cdots\subseteq I_X,
$$
where $\widetilde{K}_i(z,I_X)=\{f\in I_X:\deg_{x_0}f\leq i\}$, and we obtain homogeneous ideals
$$
K_i(z,I_X)=\{\partial^i f/\partial x_0^i:f\in\widetilde{K}_i(z,I_X)\}\subseteq S_z,
$$
called the {\it $i$-th partial elimination ideal}, together with an ascending chain
$$
K_0(z,I_X)\subseteq K_1(z,I_X)\subseteq\cdots\subseteq K_{s-1}(z,I_X)\subsetneq K_s(z,I_X)=K_{s+1}(z,I_X)=\cdots.
$$
The chain $K_i(z,I_X)$ eventually stabilizes at $s(z,I_X):=s$ to a homogeneous ideal $K_\infty(z,I_X):=K_s(z,I_X)$ in $S_z$.

Then, $K_0(z,I_X)$ is the homogeneous ideal of $X_z\subseteq\mathbb P^{r-1}$, and $K_\infty(z,I_X)$ defines $\mathbb PTC_zX$ scheme-theoretically. Moreover, by Proposition \ref{geodescription} and Corollary \ref{sprime} in the appendix, we have the set-theoretical equality
$$
V(K_i(z,I_X))=\{p\in\mathbb P^{r-1}:\textup{length}(\pi_z^{-1}(p))>i\}\phantom{.}\cup\phantom{.}\mathbb PTC_zX,
$$
where $\pi_z$ is the projection map $X\setminus\{z\}\to X_z$, and if $\pi_z$ is generically finite, then
$$
\deg \pi_z=\min\{i:\dim V(K_i(z,I_X))<\dim X_z\}.
$$
For details, refer to the appendix.

\bigskip

\noindent{\bf Tangential projections.} The projection $\pi_{\mathbb T_zX}:X\dashrightarrow X_{\mathbb T_zX}\subseteq\mathbb P^{r-n-1}$, $n:=\dim\mathbb T_zX$, of $X$ from a projective tangent space $\mathbb T_zX$ is called a {\it tangential projection}, and repeating tangential projections $m$ times, we obtain an {\it $m$-tangential projection}. A general $m$-tangential projection $X_\Lambda$ will often be used, where $\Lambda=\langle\mathbb T_{z_1}X,\ldots,\mathbb T_{z_m}X\rangle$ for $m$ general points $z_1,\ldots,z_m$ of $X$.

Recall the well-known {\it Terracini's lemma}:
$$
\mathbb T_wS^q(X)=\langle\mathbb T_{z_1}X,\ldots,\mathbb T_{z_q}X\rangle
$$
for $q$ general points $z_1,\ldots,z_q\in X$ and a general point $w\in\langle z_1,\ldots,z_q\rangle\subseteq S^q(X)$. Let $\Lambda$ be the span of $m$ general projective tangent spaces to $X$ with $\lambda:=\dim\Lambda<r$ and $m<q$. Consider the $m$-tangential projection $X_\Lambda\subset\mathbb P^{r-\lambda-1}$ and its $(q-m)$-secant variety $S^{q-m}(X_\Lambda)$. By Terracini's lemma, $S^{q-m}(X_\Lambda)\subseteq\mathbb P^{r-\lambda-1}$ has the same codimension as $S^q(X)\subseteq\mathbb P^r$ (cf.\phantom{.}\cite[Lemma 1.13]{ciliberto2006varieties}).

 The general $m$-tangential projection is a basic ingredient of both the tangent cone structure of higher secant varieties and the generalized Bronowski's conjecture.

\bigskip

\noindent{\bf The structure of tangent cones.} Note that whenever $q\geq 2$ and $S^q(X)\neq\mathbb P^r$, $S^q(X)$ is singular at every point of $S^m(X)$, $1\leq m<q$. Understanding the tangent cone to $S^q(X)$ at a general point $z\in S^m(X)$ is very important in geometry and syzygies of higher secant varieties. The structure of tangent cone in connection with the tangential projection is the origin of matryoshka structure of higher secant varieties in many aspects.

\begin{Prop}[Ciliberto-Russo, {\cite[Theorem 3.1]{ciliberto2006varieties}}]\label{tancone}
Let $1\leq m<q$, and $z_1,\ldots,z_m$ $m$ general points of $X$. Choose a general point $z$ in $\langle z_1,\ldots,z_m\rangle\subseteq S^m(X)$, and write $\Lambda=\langle\mathbb T_{z_1}X,\ldots,\mathbb T_{z_m}X\rangle$. Then, the cone over $S^{q-m}(X_\Lambda)$ with vertex $\Lambda$ is an (reduced) irreducible component of the projective tangent cone $\mathbb TC_zS^q(X)$.
\end{Prop}

\begin{Def}
One can say that $S^q(X)$ {\it has simple tangent cone} at $z\in X$ if $s(z,I_{S^q(X)})=1$ and $K_\infty(z,I_{S^q(X)})$ is prime.
\end{Def}

\begin{Rmk}
Let $z$ be a general point in $X$, and suppose that $S^q(X)$ has simple tangent cone at $z$. If $q=1$, then $K_1(z,I_X)$ is generated by linear forms which define $\mathbb PT_zX$. For $q\geq 2$, similar observation describes the following:
\begin{enumerate}[\normalfont (1)]
\item  By the partial elimination ideal theory and Proposition \ref{tancone}, we obtain a chain of homogeneous ideals
$$
K_1(z,I_{S^q(X)})\subseteq\cdots\subseteq K_\infty(z,I_{S^q(X)})\subseteq I_{\mathbb PTC_zS^q(X)}\subseteq(I_{S^{q-1}(X_{\mathbb T_zX})}),
$$
where $I_{\mathbb PTC_zS^q(X)}$ is the (saturated) homogeneous ideal defining the projective scheme $\mathbb PTC_zS^q(X)$ and $(I_{S^{q-1}(X_{\mathbb T_zX})})\subseteq S_z$ is the prime ideal generated by $I_{S^{q-1}(X_{\mathbb T_zX})}$.

\item If $S^q(X)$ has simple tangent cone at the general point $z\in X$, then all the homogeneous ideals in the chain above are the same so that $\mathbb TC_zS^q(X)$ is irreducible and reduced, and in fact, equal to the cone over $S^{q-1}(X_{\mathbb T_zX})$, and moreover, the projection $\pi_z:S^q(X)\dashrightarrow S^q(X)_z=S^q(X_z)$ is birational (Corollary \ref{sprime}).

\item It is the case for higher secant varieties of minimal degree and del Pezzo higher secant varieties; see Proposition \ref{paring} and Proposition \ref{paringII}.
\end{enumerate}
\end{Rmk}

\begin{Ex}
There are non-examples of having simple tangent cone.

(1) Let $C$ be a trigonal smooth curve of genus $3$. Take $K$ to be the canonical divisor on $C$, and $D$ the divisor associated to the trigonal map $C\to\mathbb P^1$. Embed $C\subset\mathbb P^4$ by $|K+D|$. Macaulay2 presents that the projectivized tangent cone $\mathbb PTC_zS^2(C)$ to a general point $z\in C$ is the union of the cone over $C_{\mathbb T_zC}$ and a linear space. For details, when the general point $z\in C$ has coordinates $(x_0:x_1:\cdots:x_4)=(1:0:\cdots:0)$, $S^2(C)$ is a hypersurface $f=0$ of degree $12$ with $\deg_{x_0}f=6$. Therefore, $\deg\{\pi_z:S^2(C)\to S^2(C_z)=\mathbb P^3\}=s(z,I_{S^2(C)})=6$. And we have $\deg C_{\mathbb T_zC}=5$. Similarly, using a trigonal smooth curve of genus $4$, we have another non-example in $\mathbb P^5$.

(2) Let $S\subset\mathbb P^8$ be a rational normal surface (for instance, $S(3,4)$). Take $X\subset\mathbb P^7$ to be the projection image of $S$ from a general point in $S^2(S)$. Then, by computing via Macaulay2 (\cite{M2}), for a general point $z\in X$,
$$
K_0(z,I_{S^2(X)})\subsetneq K_1(z,I_{S^2(X)})\subsetneq K_2(z,I_{S^2(X)})=K_\infty(z,I_{S^2(X)})=(I_{X_{\mathbb T_zX}}).
$$
Hence, $K_\infty(z,I_{S^2(X)})$ is prime, but $s(z,I_{S^2(X)})=2$. For more information, see Proposition \ref{almost}.
\end{Ex}

\bigskip

From the degree formula \eqref{degformula} and the fact that we obtain $\deg\mathbb TC_zS^q(X)\geq\deg S^{q-1}(X_{\mathbb T_zX})$ for a general point $z\in X$ (Proposition \ref{tancone}), the results below follow:

\begin{Prop}[cf.\phantom{.}{\cite[Lemma 4.1]{ciliberto2006varieties}}]
Let $z\in X$ be a general point. Then, we have
$$
\deg S^q(X)\geq\deg S^q(X_z)+\deg S^{q-1}(X_{\mathbb T_zX}),
$$
with equality if and only if $\deg \mathbb TC_zS^q(X)=\deg S^{q-1}(X_{\mathbb T_zX})$ with the projection $\pi_z:S^q(X)\dashrightarrow S^q(X_z)$ birational.
\end{Prop}

\begin{Cor}[{\cite[Theorem 4.2]{ciliberto2006varieties}}]\label{degbound}
We have
$$
\deg S^q(X)\geq\binom{e+q}{q}.
$$
And the equality implies that $\deg S^q(X_z)=\binom{e-1+q}{q}$ and $\deg S^{q-1}(X_{\mathbb T_zX})=\binom{e+q-1}{q-1}$ for a general point $z$ in $X$.
\end{Cor}

\begin{Def}
One says that $S^q(X)$ is a {\it $q$-secant variety of minimal degree} if $\deg S^q(X)$ has its own minimal value, i.e., $\deg S^q(X)=\binom{e+q}{q}$ for $e:=\text{codim}S^q(X)$. And we call $X$ an {\it $\textup{M}^q$-variety} if $S^q(X)$ has minimal degree, and an {\it $\textup{MA}^q$-variety} if it is an $\text{M}^q$-variety and the map $\widetilde{S}^q(X)\to S^q(X)$ is birational.
\end{Def}

We use slightly different notations from those in \cite{ciliberto2006varieties}. In the definition above, we do not exclude the case $S^q(X)=\mathbb P^r$. 

\begin{Rmk}
The projective variety $X$ need not be an $\text{MA}^q$-variety even if it is an $\text{M}^q$-variety and $\dim S^q(X)=\dim\widetilde{S}^q(X)<\dim\mathbb P^r$. Indeed, take $X$ to be the image of $\mathbb P^2$ in $\mathbb P^6$ induced by the linear system $|6L-2p_1-\cdots-2p_7|$, where $L$ is any line in $\mathbb P^2$ and $p_1,\ldots,p_7$ are seven general points of $\mathbb P^2$. Then, by using \cite[Theorem 9.1]{ciliberto2006varieties}, $X$ is an $\text{M}^2$-variety with $\dim S^2(X)=\dim\widetilde{S}^2(X)$. On the other hand, let $p_8$ be a general point of $\mathbb P^2$ with $z_8$ its image in $X$. Then, $X_{z_8}$ is the image of $\mathbb P^2$ in $\mathbb P^5$ corresponding to the linear system $|6L-2p_1-\cdots-2p_7-p_8|$. By \cite[Theorem 1.1]{ciliberto2005surfaces}, $X_{z_8}$ is not an $\text{MA}^2$-variety, hence neither is $X$.
\end{Rmk}

\begin{Rmk}
Varieties with one apparent double point, i.e., $\text{MA}^2$-varieties $X^n\subset\mathbb P^{2n+1}$, are a classical object, and have been studied for a long time. Note that {\it smooth} $\text{MA}^2$-varieties $X^n\subset\mathbb P^{2n+1}$ of dimension $n$ are classified and listed for $n=2,3$ in \cite[Theorem 4.10 and Theorem 7.1]{ciliberto2004varieties}, and for the case $\deg X\leq 2n+4$ in \cite[Theorem 2 and Theorem 3]{alzati2003special}. 
\end{Rmk}

Closing this subsection, we summarize facts used frequently and implicitly as follows:
\begin{itemize}
\item For two general points $z,w\in X$, inner projection and tangential projection procedures commute:
\begin{enumerate}[\normalfont (1)]
\item $(X_z)_{\overline{w}}=(X_w)_{\overline{z}}=X_{\langle z,w\rangle}$;

\item $(X_z)_{\mathbb T_{\overline{w}}(X_z)}=(X_{\mathbb T_wX})_{\overline{\overline{z}}}=X_{\langle z,\mathbb T_wX\rangle}$; and

\item $(X_{\mathbb T_zX})_{\mathbb T_{\overline{\overline{w}}}(X_{\mathbb T_zX})}=(X_{\mathbb T_wX})_{\mathbb T_{\overline{\overline{z}}}(X_{\mathbb T_wX})}=X_{\langle\mathbb T_zX,\mathbb T_wX\rangle}$,
\end{enumerate}
where $\overline{z}=\pi_w(z)$, $\overline{\overline{z}}=\pi_{\mathbb T_wX}(z)$, $\overline{w}=\pi_z(w)$ and $\overline{\overline{w}}=\pi_{\mathbb T_zX}(w)$.
$$
\xymatrix@C=2cm{
X_{\langle z,\mathbb T_wX\rangle} & X_{\mathbb T_wX} \ar@{-->}[l]_--{\text{inn.}} \\
X_z \ar@{-->}[u]^--{\text{tan.}} & X \ar@{-->}[l]^--{\text{inn.}} \ar@{-->}[u]_--{\text{tan.}} \ar@{-->}[lu]
}
\quad
\xymatrix@C=2cm{
X_{\langle \mathbb T_zX,\mathbb T_wX\rangle} & X_{\mathbb T_wX} \ar@{-->}[l]_--{\text{tan.}} \\
X_{\mathbb T_zX} \ar@{-->}[u]^--{\text{tan.}} & X \ar@{-->}[l]^--{\text{tan.}} \ar@{-->}[u]_--{\text{tan.}} \ar@{-->}[ul]
}
$$

\item For any linear space $\Lambda\neq\mathbb P^r$, we have $S^q(X)_{\Lambda}=S^q(X_\Lambda)$.

\item If $e=\text{codim}S^q(X)$, then for a general point $z$ in $X$, $S^q(X_z)$ has codimension $e-1$ and $S^{q-1}(X_{\mathbb T_zX})$ has codimension $e$.
\end{itemize}

\subsection{Equations and syzygies}

In this subsection, we will collect special properties on equations and syzygies of higher secant varieties, which in general, do not hold in the category of all projective varieties. To the end, let us introduce some definitions: Let $M$ be a finitely generated graded $S$-module, and consider the minimal free resolution
$$
\resizebox{\textwidth}{!}{
\xymatrix{
0 & M \ar[l] & F_0 \ar[l] & F_1 \ar[l] & \cdots \ar[l] & F_i=\bigoplus_jS^{\beta_{i,j}(M)}(-i-j) \ar[l] & \cdots. \ar[l]
}
}
$$
Here, the $\beta_{i,j}(M)$ are called the {\it (graded) Betti numbers} of $M$. The {\it Betti table} of $M$ is a diagram of the following form:
$$
\begin{array}{c|cccc}
               & \cdots & i                           & i+1                         & \cdots \\ \hline
\vdots &  & \vdots                  & \vdots                       &  \\
j       & \cdots & \beta_{i,j}(M)     & \beta_{i+1,j}(M)     & \cdots \\
j+1     & \cdots & \beta_{i,j+1}(M) & \beta_{i+1,j+1}(M) & \cdots \\
\vdots &  & \vdots                   & \vdots                      &
\end{array}
$$
When $\beta_{i,j}(M)=0$, we put ``$-$" at the place in the Betti table. The Betti table of $X\subseteq\mathbb P^r$ is that of the homogeneous coordinate ring $S_X$.

\begin{Def}
The homogeneous coordinate ring $S_X$ is said to {\it satisfy property $N_{d,p}$} if $\beta_{i,j}(S_X)=0$ for every $i\leq p$ and every $j\geq d$. Also,$X$ is {\it $m$-regular} if $\beta_{i,j}(S_X)=0$ for every $j\geq m$, and we define the {\it regularity} of $X$ as $\text{reg}X:=\min\{m:\text{$X$ is $m$-regular}\}$.
\end{Def}

\begin{Def}
Let $M$ be a graded $S$-module. We define the {\it Koszul cohomology group $K_{i,j}(M,V)$} to be the cohomology of the complex
\begin{displaymath}
\xymatrix{
\bigwedge^{i+1}V\otimes M_{j-1} \ar[r]^--{\delta_{i+1,j-1}} & \bigwedge^iV\otimes M_j \ar[r]^--{\delta_{i,j}} & \bigwedge^{i-1}V\otimes M_{j+1}
}
\end{displaymath}
at the middle, where the differential $\delta_{i,j}$ is defined by
$$
\delta_{i,j}(v_1\wedge\cdots\wedge v_i\otimes m)=\sum_{k=1}^i(-1)^{k-1}v_1\wedge\cdots\wedge\widehat{v_k}\wedge\cdots\wedge v_i\otimes v_km.
$$
An element of $K_{i,j}(M,V)$ is called a {\it Koszul class}.
\end{Def}

\bigskip

\noindent{\bf Conventions.}
We introduce several conventions:
\begin{itemize}
\item $\beta_{i,j}(X)=\beta_{i,j}(S_X)$.

\item $K_{i,j}(M)=K_{i,j}(M,V)$ if $V$ is well understood.

\item $K_{i,j}(X)=K_{i,j}(S_X)$.

\item $I_\emptyset$ is the irrelevant ideal $S_+$ of $S$ for the empty set $\emptyset\subset\mathbb P^r$. (Recall that $S^0(X)=\emptyset\subset\mathbb P^r$.)

\item $S_\emptyset$ is the residue field $S/I_\emptyset$ of $S$.

\item $\beta_{i,j}(\emptyset)=\beta_{i,j}(S_\emptyset)$ and $K_{i,j}(\emptyset)=K_{i,j}(S_\emptyset)$.

\end{itemize}
Note that the residue field $S_\emptyset$ is minimally resolved by the Koszul complex
$$
\xymatrix{
0 & S_\emptyset \ar[l] & S \ar[l] & V\otimes S(-1) \ar[l] & \bigwedge^2V\otimes S(-2) \ar[l] & \cdots. \ar[l]
}
$$

\bigskip

\noindent{\bf Equations and syzygies of $S^q(X)$.} Now, some observations on equations and syzygies of higher secant varieties are made as follows:

\begin{Lem}\label{noqform}
Take linearly independent $r+1$ points $z_0,\ldots,z_r$ in $X\subseteq\mathbb P^r$, and their dual basis $x_0,\ldots,x_r$ for $V$. Then, every monomial of an equation for $S^q(X)$ consists of more than $q$ distinct variables $x_i$. In particular, no hypersurfaces of degree $q$ contain $S^q(X)$, and we have $\deg_{x_0}f\leq m-q$ for all $f\in(I_{S^q(X)})_m$.
\end{Lem}

\begin{proof}
Since an equation of $S^q(X)$ always vanishes on linear spans of linearly independent $q$ points in $X$, we are done.
\end{proof}

Therefore, the Betti table of $S^q(X)$ can not have a nonzero entry in the $i$-th row for any $i\leq q-1$ except $\beta_{0,0}$:
$$
\begin{array}{c|cccccc}
& 0 & 1 & 2 & \cdots & p & \cdots \\ \hline
0 & 1 & - & - & \cdots & - & \cdots \\
1 & - & - & - & \cdots & - & \cdots \\
\vdots & \vdots & \vdots & \vdots & & \vdots & \\
q-1 & -  & - & - & \cdots & - & \cdots \\
q & - & \beta_{1,q} & \beta_{2,q} & \cdots & \beta_{p,q} & \cdots \\
\vdots & \vdots & \vdots & \vdots & & \vdots &
\end{array}
$$

We have the following basic inequality on the Betti numbers appearing in each first possible nontrivial row of higher secant varieties, which is a consequence of Proposition \ref{basicA}, Lemma \ref{noqform} and the containment $K_1(z,I_{S^q(X)})\subseteq(I_{S^{q-1}(X_{\mathbb T_zX})})$ with no degree $q-1$ equations of $S^{q-1}(X_{\mathbb T_zX})$:

\begin{Prop}[The basic inequality]\label{basic}
Taking a general point $z\in X$, we have
$$
\beta_{p,q}(S^q(X))\leq\beta_{p,q}(S^q(X_z))+\beta_{p-1,q}(S^q(X_z))+\beta_{p,q-1}(S^{q-1}(X_{\mathbb T_zX})).
$$
And if $S^q(X)$ has simple tangent cone at a general point of $X$, and if $\beta_{p-1,q+1}(S^q(X_z))$ and $\beta_{p-1,q}(S^{q-1}(X_{\mathbb T_zX}))$ vanish, then the equality holds.
\end{Prop}

For the case $q=1$, refer to \cite[Theorem 3.1]{han2015sharp}. Notice that according to our conventions, if $q=1$, then the value of $\beta_{p,q-1}(S^{q-1}(X_{\mathbb T_zX}))$ is equal to $\binom{e}{p}$, $e=\text{codim}X$.

\begin{Cor}[Upper bounds of Betti numbers]\label{syzbound}
Let $e=\textup{codim}S^q(X)$. Then, we have
$$
\beta_{p,q}(S^q(X))\leq B^e_{p,q}:=\binom{p+q-1}{q}\binom{e+q}{p+q}
$$
for every $p\geq 1$. And if the equality holds for some $p\geq 1$, then for a general point $z$ in $X$, we have $\beta_{p,q}(S^q(X_z))=B^{e-1}_{p,q}$, $\beta_{p-1,q}(S^q(X_z))=B^{e-1}_{p-1,q}$ and $\beta_{p,q-1}(S^{q-1}(X_{\mathbb T_zX}))=B^e_{p,q-1}$.
\end{Cor}

\begin{proof}
We will use double induction on $e$ and $q$. The case $q=1$ follows from \cite[Theorem 1.2]{han2015sharp}, and the case $e=1$ is trivial. Due to the equality $B^e_{p,q}=B^{e-1}_{p,q}+B^{e-1}_{p-1,q}+B^e_{p,q-1}$, we are done.
\end{proof}

In addition, we present a useful theorem on Koszul cohomology vanishing. This is very important in proving the syzygetic characterization theorems: Theorems \ref{syzM} and \ref{syzdP}. For its proof, refer to Proposition \ref{vanishingA}.

\begin{Thm}[Vanishing of Betti numbers]\label{vanishing}
Suppose that $S^q(X)$ has simple tangent cone at a general point $z$ of $X$. If
\begin{enumerate}[\normalfont (1)]
\item $\beta_{i,j-1}(S^{q-1}(X_{\mathbb T_zX}))=0$; and

\item $\beta_{i-1,j}(S^q(X_z))=\beta_{i,j}(S^q(X_z))=0$,
\end{enumerate}
then we have $\beta_{i,j}(S^q(X))=0$.
\end{Thm}

\section{The generalized Bronowski's conjecture}

As explained in the introduction, the generalized Bronowski's conjecture is raised by C. Ciliberto and F. Russo based on Bronowski's claim, and reads off algebraic and geometric properties of a $q$-secant variety from those of a general $(q-1)$-tangential projection. Recall that the conjecture asserts the following: If a general $(q-1)$-tangential projection of $X$ is birational, and if its image is a variety of minimal degree, then $X$ is an $\textup{MA}^q$-variety. The converse was already proved in \cite[Corollay 4.5]{ciliberto2006varieties}.

\subsection{A proof of the generalized Bronowski's conjecture, weak form}

\begin{Def}
For simplicity, we define the following:
\begin{enumerate}[\normalfont (1)]
\item $S^q(X)$ satisfies {\it property $ST$} if it has simple tangent cone at a general point of $X$.

\item $S^q(X)$ satisfies {\it property $IC$} if
$
I_{S^q(X)}=\sum_{i=0}^{q+1}(I_{S^q(X_{z_i})})
$
for $q+2$ general points $z_i$ in $X$.
\end{enumerate}
\end{Def}

\begin{Ex}
Let $X=C$ be a linearly normal smooth curve of genus $2$. Suppose that $C$ is a space curve. Then, $C$ has a (unique) quadratic equation, but neither does any general inner projection of $C$. Thus, $C$ does not satisfy property $IC$. Now, suppose that $\text{codim}C=3$. It is known that $I_C$ is minimally generated by four quadratic forms. Since a general inner projection of $C$ lies in only one quadric, the sum $(I_{X_{z_0}})+(I_{X_{z_1}})+(I_{X_{z_2}})$, $z_i\in X$ general, has at most three quadratic forms, and thus $C$ does not satisfies property $IC$. If $\text{codim}C\geq 4$, then $C$ satisfies property $IC$ by Corollary \ref{paringCor}.
\end{Ex}

The following proposition shows that property $ST$ and property $IC$ are very closely related, and is a key step towards proving the generalized Bronowski's conjecture, weak form:

\begin{Prop}\label{paring}
Let $e=\textup{codim}S^q(X)$. The following hold:
\begin{enumerate}[\normalfont (1)]
\item Set $e\geq 2$ and $q=1$. Suppose that $X_z$ satisfies property $ST$ for a general point $z$ in $X$. Then, $X$ satisfies both property $ST$ and property $IC$.

\item Set $e\geq 2$ and $q\geq 2$. Suppose that for a general point $z$ in $X$, $S^q(X_z)$ satisfies property $ST$ and $S^{q-1}(X_{\mathbb T_zX})$ satisfies property $IC$. Then, $S^q(X)$ satisfies both property $ST$ and property $IC$.

\item Set $e\geq 1$ and $q\geq 1$. If $S^q(X)$ is a $q$-secant variety of minimal degree, then $S^q(X)$ satisfies property $ST$, and if in addition, $e\geq 2$, then $S^q(X)$ satisfies property $IC$.
\end{enumerate}
\end{Prop}

\begin{proof}
In the argument below, $\overline{z}$ denotes the image of $z$ under an appropriate projection map, and $\mathbb T_i$ projective tangent spaces to $X$ at $z_i\in X$ when an indexed set $z_i$ of points is given.

(1) A proof of (1) is similar to that of (2) below because for three general points $z_0,z_1,z_2$ in $X$, we obtain $\mathbb T_0=\langle z_1,\mathbb T_0\rangle\cap\langle z_2,\mathbb T_0\rangle$, and so
$$
V_{\mathbb T_0}=V_{\langle z_1,\mathbb T_0\rangle}+V_{\langle z_2,\mathbb T_0\rangle},
$$
where $V_\Lambda$ is the vector space of linear forms vanishing on a linear space $\Lambda$.

(2) Set $q\geq 2$. Take $q+2$ general points $z_i$ in $X$, and say $z_0=(1:0:\cdots:0)$ in homogeneous coordinates $(x_0:x_1:\cdots:x_r)$. For the property $ST$ part, choose an arbitrary element $f\in I_{S^{q-1}(X_{\mathbb T_0})}$. Note that by the assumption,
$$
f\in I_{S^{q-1}(X_{\mathbb T_0})}=\sum_{i=1}^{q+1}(I_{S^{q-1}(X_{\langle\mathbb T_0,z_i\rangle})})=\sum_{i=1}^{q+1}(K_1(\overline{z_0},I_{S^q(X_{z_i})})).
$$
So, $I_{S^q(X)}$ has an element $g$ of form $fx_0+\text{(terms without $x_0$)}$, and thus $f\in K_1(z_0,I_{S^q(X)})$. For the property $IC$ assertion, let $f\in I_{S^q(X)}$. Now, we will use induction on $d:=\deg_{x_0}f$. If $d=0$, then $f\in (I_{S^q(X_{z_0})})$. Suppose that $f=f'x_0^d+\text{(lower $x_0$-degree terms)}$ with $d\geq 1$. Then, by the assumption, $f'$ lies in $K_d(z_0,I_{S^q(X)})$, hence
$$
f'\in(I_{S^{q-1}(X_{\mathbb T_0})})=\sum_{i=1}^{q+1}(I_{S^{q-1}(X_{\langle\mathbb T_0,z_i\rangle})})=\sum_{i=1}^{q+1}(K_d(\overline{z_0},I_{S^q(X_{z_i})})).
$$
Therefore, an element $g$ of form $f'x_0^d+\text{(lower $x_0$-degree terms)}$ exists in $\sum_{i=1}^{q+1}(I_{S^q(X_{z_i})})$. Consider the difference $f-g$ and its degree in $x_0$ so that the value of $d$ decreases. Thus, $f\in\sum_{i=0}^{q+1}(I_{S^q(X_{z_i})})$. We are done.

(3) Note that if $e=1$ and $q\geq 1$, then property $ST$ is obvious since $S^q(X)$ and $S^{q-1}(X_{\mathbb T_zX})$ are both hypersurfaces, and that the case $e\geq 2,q=1$ can be proved by (1). Using (2) and double induction on $e,q\geq 2$, we are done.
\end{proof}

\begin{Cor}\label{paringCor}
If $X$ satisfies property $N_{2,2}$ with $\textup{codim}X\geq 2$, then it satisfies property $IC$.
\end{Cor}

\begin{proof}
By (1) of the proposition above, it is enough to show that a general inner projection $X_z$ satisfies property $ST$. By \cite{han2012analysis}, $X_z$ satisfies property $N_{2,1}$, hence property $ST$.
\end{proof}

The following lemma is useful to deal with the case where $S^q(X)$ has small codimension. In particular, it is used for the proof of Theorem \ref{gBcwf}.

\begin{Lem}\label{prolongation}
Let $q\geq 2$. Given an indexed set $z_i$ of points in $X$, set a homogeneous coordinate system $x_0,\ldots,x_r$ such that $V(x_j)\ni z_i$ if $i\neq j$, and put $\mathbb T_i:=\mathbb T_{z_i}X$. Suppose that
\begin{enumerate}[\normalfont (1)]
\item every $(q-1)$-form vanishing on $S^{q-2}(X_{\langle \mathbb T_0,\mathbb T_1\rangle})$ is of the form $\partial f/\partial x_1$ with some $q$-form $f$ vanishing on $S^{q-1}(X_{\mathbb T_0})$ for two general points $z_0,z_1$ in $X$;
and
\item there are no degree $q-1$ hypersurfaces containing $S^{q-2}(X_{\langle z_0,\mathbb T_1,\mathbb T_2\rangle})$ for three general points $z_0,z_1,z_2$ in $X$.
\end{enumerate}
Then, every $q$-form vanishing on $S^{q-1}(X_{\mathbb T_0})$ is of the form $\partial f/\partial x_0$ with some $(q+1)$-form $f$ vanishing on $S^q(X)$ for a general point $z_0$ in $X$.
\end{Lem}

\begin{proof}
Choose $r+1$ general points $z_i$ in $X$, and let $f_0$ be any $q$-form in $I_{S^{q-1}(X_{\mathbb T_0})}$. Then, $\deg_{x_i}f_0\leq 1$ and
$$
\frac{\partial f_0}{\partial x_i}\in (I_{\mathbb PTC_{\overline{z_i}}S^{q-1}(X_{\mathbb T_0})})_{q-1}\subseteq(I_{S^{q-2}(X_{\langle \mathbb T_0,\mathbb T_i\rangle})})_{q-1},
$$
which by assumption, says that for each $1\leq i\leq r$, there is some $q$-form $f_i\in I_{S^{q-1}(X_{\mathbb T_i})}$ such that $\partial f_i/\partial x_0=\partial f_0/\partial x_i$.

We claim that
$$
\frac{\partial f_i}{\partial x_j}=\frac{\partial f_j}{\partial x_i}
$$
for every pair $i\neq j$. Indeed, we have
$$
\frac{\partial^2f_i}{\partial x_0\partial x_j}=\frac{\partial^2 f_i}{\partial x_j\partial x_0}=\frac{\partial^2 f_0}{\partial x_j\partial x_i}=\frac{\partial^2 f_0}{\partial x_i\partial x_j}=\frac{\partial^2 f_j}{\partial x_i\partial x_0}=\frac{\partial^2 f_j}{\partial x_0\partial x_i},
$$
and by assumption in consideration of the difference $\partial f_i/\partial x_j-\partial f_j/\partial x_i$, the claim is verified. Now, by the claim above, the $(q+1)$-form
$$
f:=\frac{1}{q+1}\sum_{i=0}^rx_if_i
$$
satisfies $\partial f/\partial x_i=f_i\in I_{S^{q-1}(X_{\mathbb T_i})}\subset I_{S^{q-1}(X)}$. Due to the prolongation theorem \cite[Theorem 1.2]{sidman2009prolongations}, we obtain $f\in I_{S^q(X)}$ with $\partial f/\partial x_0=f_0$. We are done.
\end{proof}

We are ready to prove the generalized Bronowski's conjecture, weak form.

\begin{proof}[\bf Proof of Theorem \ref{gBcwf}]
(1) We use induction on $q\geq 1$. For $q=1$, it is trivial that if $X_z$ is a variety of minimal degree, then so is $X$. Now, suppose that $q\geq 2$. Let $w\in X$ be another general point. Then, by the induction hypothesis, $S^{q-1}(X_{\mathbb T_zX})$ has minimal degree since so does $S^{q-1}(X_{\langle w,\mathbb T_zX\rangle})$. By Proposition \ref{paring}, $S^q(X)$ has simple tangent cone at a general point of $X$. Hence, $S^q(X)$ is a $q$-secant variety of minimal degree.

(2) Similarly, use induction on $e\geq 1$. The case $e=1$ is due to Lemma \ref{prolongation}. Let $z,w\in X$ be general points, and consider $S^{q-1}(X_{\mathbb T_wX})$, $S^{q-1}(X_{\langle z,\mathbb T_wX\rangle})$ and $S^q(X_z)$. As above, $S^q(X_z)$ has minimal degree as well as $S^{q-1}(X_{\mathbb T_zX})$, and $S^q(X)$ has simple tangent cone at a general point in $X$. We are done.
\end{proof}

\begin{Cor}\label{gBcwfCor}
A projective variety $X$ is an $\textup{M}^q$-variety if and only if its general $(q-1)$-tangential projection is a minimal degree variety.
\end{Cor}

\begin{proof}
The ``only if'' part is by Corollary \ref{degbound}. In order to show the ``if'' part, we use induction on $q$. The case $q=1$ is trivial. Now, let $q\geq 2$, and $X'$ a general tangential projection of $X$. Then, by assumption, a general $(q-2)$-tangential projection of $X'$ is a variety of minimal degree. The induction hypothesis says that $X'$ is an $M^{q-1}$-variety. We are done by Theorem \ref{gBcwf}.
\end{proof}

As a general tangential projection of a minimal degree variety again has minimal degree, Corollary \ref{gBcwfCor} implies the following:

\begin{Cor}\label{minmin}
If $X$ is an $\textup{M}^q$-variety, then it is also an $\textup{M}^{q+1}$-variety. In particular, any variety of minimal degree is an $\textup{M}^q$-variety for all $q\geq 1$.
\end{Cor}

The rest of this subsection will be devoted to revealing the portion of the weak form in solving the generalized Bronowski's conjecture.

\begin{Rmk}
The generalized Bronowski's  conjecture is true for the following cases:
\begin{enumerate}[\normalfont (1)]
\item $e=0,q=2$ and $X$ is smooth with $1\leq\dim X\leq 3$; and
\item $e\geq0,q\geq 2$ and $X$ is a smooth surface,
\end{enumerate}
where $e=\text{codim}S^q(X)$. See \cite[pp.\phantom{.}26-27]{ciliberto2006varieties}. Their arguments depend on the classification of $\text{MA}^q$-varieties.
\end{Rmk}

\begin{Cor}\label{e=0}
If the generalized Bronowski's conjecture holds for the case $S^q(X)=\mathbb P^r$, then so does it for all the cases.
\end{Cor}

\begin{proof}
Let $e\geq 1$ be the codimension of $S^q(X)$. Suppose that a general $(q-1)$-tangential projection of $X$ is birational, and its image $X_\Lambda$ is a variety of minimal degree. Let $X'$ be a general $e$-inner projection of $X$. Since a general $e$-inner projection of $X_\Lambda$ is birational, and since two projections commute, a general $(q-1)$-tangential projection of $X'$ is birational onto the projective space. Now, if the generalized Bronowski's conjecture holds for $X'$, then $X'$ is an $\text{MA}^q$-variety, hence $X$ is an $\text{MA}^q$-variety by Theorem \ref{gBcwf} and the commutative diagram
$$
\xymatrix@C=2cm{
\widetilde{S}^q(X) \ar@{-->}[r]^{\substack{\text{generically} \\ 1:1}} \ar[d] & \widetilde{S}^q(X') \ar[d]^{\substack{\text{generically} \\ 1:1}} \\
S^q(X) \ar@{-->}[r]_--{\substack{\text{generically} \\ 1:1}} & S^q(X')
}
$$
(cf. \cite[Lemma 4.1 (vii)]{ciliberto2006varieties}).
\end{proof}

Now, we will prove the generalized Bronowski's conjecture for curves. The following definition will be used:

\begin{Def}[\cite{kaji1986tangentially}]
Let $X=C\subset\mathbb P^r$ be a curve. One says that $C$ is {\it tangentially degenerate} if for a general point $z\in C$, there is a point $z'\in C\cap\mathbb T_zC$ other than $z$.
\end{Def}

\begin{Cor}\label{gBccurve}
The generalized Bronowski's conjecture is true for curves.
\end{Cor}

\begin{proof}
Let $X=C\subseteq\mathbb P^r$ be a curve. If $S^q(C)\neq\mathbb P^r$, then this is because of Theorem \ref{gBcwf} and the fact that the map $\widetilde{S}^q(C)\to S^q(C)$ is birational with $S^q(C)\neq\mathbb P^r$ (cf.\phantom{.}\cite[Proposition 1.5]{ciliberto2006varieties}). We only need to consider the curve $C$ in $\mathbb P^{2q-1}$ such that its general $(q-1)$-tangential projection is birational onto $\mathbb P^1$.  Let $C'\subset\mathbb P^3$ be a general $(q-2)$-tangential projection of $C$. Write $d=\deg C'$. As a general tangential projection of $C'$ is birational, we have $\text{length}(C'\cap\mathbb T_zC')=d-1$ for a general point $z\in C'$, which implies that $C'$ has maximal regularity, namely $\text{reg}C'=d-1$. Now, by referring to the classification of curves with maximal regularity in \cite{gruson1983theorem}, if $d>4$, then $C'$ is smooth, and so not tangentially degenerate (\cite[Theorem 3.1]{kaji1986tangentially}), which means that $\text{length}(C'\cap\mathbb T_zC')=2$, a contradiction. Note that the following is well known: any nondegenerate irreducible projective curve of almost minimal degree is smooth or cut out by quadrics. So, by a similar reason, $C'$ can not have almost minimal degree. Thus, $C'$ is a twisted cubic curve, hence $C$ is an $\text{M}^{q-1}$-variety by Theorem \ref{gBcwf}. Thus, as noted in \cite[Theorem 6.1]{ciliberto2006varieties}, $C$ is an $\text{MA}^q$-variety since it is a rational normal curve.
\end{proof}

For a result applied to any dimensions, refer to the remark below:

\begin{Rmk}
The projective variety $X$ is said to be {\it $q$-tangentially weakly defective} if $q$ general points $z_1,\ldots,z_q$ in $X$ are not isolated points of the locus
$$
\overline{\{z\in\text{Sm}X:\mathbb T_zX\subseteq\langle\mathbb T_{z_1}X,\ldots,\mathbb T_{z_q}X\rangle\}}\subseteq X.
$$
The tangentially weak defectiveness enjoys the following properties:
\begin{enumerate}[\normalfont (1)]
\item If $X$ is not $q$-tangentially weakly defective, then it is {\it $q$-identifiable}, that is, the map $\widetilde{S}^q(X)\to S^q(X)$ is birational; see \cite[Proposition 14]{casarotti2019non}.

\item When a general $m$-tangential projection of $X$ is generically finite for an integer $1\leq m\leq q-1$, the projection image is $(q-m)$-tangentially weakly defective if and only if $X$ is $q$-tangentially weakly defective; see \cite[Proposition 14]{casarotti2020tangential}.

\item If $X$ is smooth, then it is not $1$-tangentially weakly defective; see \cite[Theorem 0]{zak1983projections}.
\end{enumerate}
\end{Rmk}

\begin{Thm}\label{gBcusual}
If $\dim S^{q+1}(X)=\dim\widetilde{S}^{q+1}(X)$, then $X$ verifies the generalized Bronowski's conjecture for $S^q(X)$. Especially, if $\dim S^m(X)=\dim\widetilde{S}^m(X)$ for some integer $m\geq q+1$, and if $X$ satisfies the assumption of generalized Bronowski's conjecture for $S^q(X)$, then $X$ is an $\textup{MA}^j$-variety for $q\leq j\leq m-1$.
\end{Thm}

\begin{proof}
Theorem \ref{gBcwf} tells us that $X$ is an $\text{M}^q$-variety. It suffices to show that the map $\widetilde{S}^q(X)\to S^q(X)$ is birational. Let $X_\Lambda$ be a general $(q-1)$-tangential projection. Then, by assumption, $X_\Lambda$ is a variety of minimal degree that satisfies $\dim S^2(X_\Lambda)=2\dim X_\Lambda+1$. This means that $X_\Lambda$ is a smooth rational normal scroll. By the remark above, $X_\Lambda$ is not $1$-tangentially defective, and thus $X$ is not $q$-tangentially weakly defective so that the map $\widetilde{S}^q(X)\to S^q(X)$ is birational.
\end{proof}

Before closing this section, it is worth mentioning the following:

\begin{Rmk}\label{gBcRmk}
In fact, the proof of theorem above shows that $X$ is an $\text{MA}^q$-variety whenever its general $(q-1)$-tangential projection  is birational onto a smooth variety of minimal degree and of positive codimension. For the complete proof of generalized Bronowski's conjecture, the remaining part is to consider the case where a general $(q-1)$-tangential projection is either a cone over a variety of minimal degree or the whole space (in this case, $\dim S^{q+1}(X)\neq\dim\widetilde{S}^{q+1}(X)$).

For instance, when $X=S(1,6)\subset\mathbb P^8$, the generalized Bronowski's conjecture holds for $S^2(X)$, but the theorem above can not be applied since $\dim S^3(X)=7$ and $\dim\widetilde{S}^3(X)=8$ are different.
\end{Rmk}

\section{A generalized $K_{p,1}$ theorem}

In this section, we will discuss the $q$-th row
$$
\beta_{1,q}(S^q(X)), \beta_{2,q}(S^q(X)),\ldots,\beta_{p,q}(S^q(X)),\ldots.
$$
of Betti table of $S^q(X)\subset\mathbb P^r$ (or the corresponding subcomplex of the minimal free resolution), namely the {\it $(q+1)$-strand} of $S^q(X)\subset\mathbb P^r$. Since $(I_{S^q(X)})_q=0$, the $(q+1)$-strand is the possibly first nontrivial row of Betti table. The {\it length of $(q+1)$-strand} is defined to be
$$
\max\{p:\beta_{p,q}\neq 0\}.
$$

The following is a generalization of a result in \cite{green1982canonical}.

\begin{Prop}\label{monomial}
Choose $r+1$ linearly independent points $z_0,\ldots,z_r\in V^\ast$ in $X$ with dual basis $x_0,\ldots,x_r\in V$, and write $x_\alpha:=x_{i_1}\wedge\cdots\wedge x_{i_{|\alpha|}}\in\bigwedge^{|\alpha|}V$ for any $\alpha=\{i_1<\cdots<i_{|\alpha|}\}\subseteq\{0,\ldots,r\}$. And let $z\in V^\ast$ be a point of $X$.
\begin{enumerate}[\normalfont (1)]
\item For any Koszul class $\eta={\displaystyle \sum_{|\alpha|=p-1}}x_\alpha\otimes f_\alpha$ in $K_{p-1,q+1}(I_{S^q(X)},V)$, we have
$$
f_\alpha=\sum_{\substack{j_0<\cdots<j_q \\ \textup{ in } \alpha^c}}a_\alpha^{j_0,\ldots,j_q}x_{j_0}\cdots x_{j_q}
$$
for some $a_\alpha^{j_0,\ldots,j_q}$ in the base field. If in addition, we fix any nonzero Koszul class $\eta$ in $K_{p-1,q+1}(I_{S^q(X)},V)$ and choose the points $z_i$ generally in $X$, then the coefficients $a_\alpha^{j_0,\ldots,j_q}$ above are all nonzero.

\item The evaluation map at $z$ extends to a map 
$$
K_{p,q}(S^q(X))\to K_{p-1,q}(S^q(X_z)).
$$
And for a given nonzero Koszul class $\gamma\in K_{p,q}(S^q(X))$, its image under the map is nonzero when we take $z$ to be general in $X$.
\end{enumerate}
\end{Prop}

\begin{proof}
Assume that a monomial $x_i\cdot M$ with $i\in\alpha$ occurs in the expression of $f_\alpha$ (with nonzero coefficient) for some monomial $M$. Then, from the fact that $\eta$ is a Koszul class, it follows that letting $\beta:=\alpha\setminus\{i\}$ for convenience, we have
$$
\sum_{j\in\beta^c}\pm x_jf_{\{j\}\cup\beta}=0.
$$
Then, the summand $\pm x_if_\alpha$ involves the monomial $x_i^2\cdot M$. So, for some $j\in\beta^c$ with $j\neq i$ and $x_j|M$, $F_{\{j\}\cup\beta}$ involves the monomial $x_i^2\cdot (M/x_j)$, for the monomial $x_i^2\cdot M$ need to be canceled out. We have established a contradiction.

 Now, in order to show that $a_\alpha^{j_0,\ldots,j_q}\neq 0$, consider the loci
$$
U_{\alpha,J}:=\{(z'_0,\ldots,z'_r):\eta(z'_{i_1}\wedge\cdots\wedge z'_{i_{p-1}})\not\equiv0\text{ on }\langle z'_{j_0},\ldots,z'_{j_q}\rangle\}\subseteq X^{r+1},
$$
for $\alpha=\{i_1<\cdots<i_{p-1}\}\subseteq\{0,\ldots,r\}$ and $J=\{j_0<\cdots<j_q\}\subseteq\{0,\ldots,r\}\setminus\alpha$.
Note that the $U_{\alpha,J}$ are open, and that
\begin{center}
$(z_0,\ldots,z_r)\in U_{\alpha,J}$ if and only if the coefficient $a_\alpha^{j_0,\ldots,j_q}$ is nonzero.
\end{center}
Furthermore, under the natural actions, $\sigma(U_{\alpha,J})=U_{\sigma(\alpha),\sigma(J)}$ holds for each element $\sigma$ of the symmetric group $\mathfrak{S}_{r+1}$, and since $\eta\neq 0$, at least one of the $U_{\alpha,J}$ is nonempty. Thus, we see that every $U_{\alpha,J}$ is nonempty, and so the general point $(z_0,\ldots,z_r)$ in $X^{r+1}$ lies in all of the $U_{\alpha,J}$, which means that each of the coefficients $a_\alpha^{j_0,\ldots,j_q}$ is nonzero.

Now, consider the mapping
$$
\eta=\sum_{|\alpha|=p-1}x_\alpha\otimes f_\alpha\mapsto\sum_{\substack{|\beta|=p-2, \\ \beta\not\ni 0}}x_\beta\otimes f_{\beta\cup\{0\}}
$$
which is an extension of the map $v_1\wedge\cdots\wedge v_{p-1}\mapsto\sum_{i=1}^{p-1}(-1)^{i-1}v_1\wedge\cdots\wedge\widehat{v_i}\wedge\cdots\wedge v_{p-1}\cdot v_i(z)$. The image is in $K_{p-2,q+1}(I_{S^{q-1}(X_z)})$, and nonzero if $\eta\neq 0$ is fixed and $z\in X$ is general (with respect to the fixed $\eta$).
\end{proof}

\begin{Rmk}
The statement (2) of the proposition above says that under the assumption, the {\it projection map} 
$$
\text{pr}_z: K_{p-1,q+1}(I_{S^q(X)},V)\to K_{p-2,q+1}(I_{S^q(X)},V_z)$$ 
(see \cite[Section 2.2.1]{aprodu2010koszul}) factors through $K_{p-2,q+1}(I_{S^q(X_z)},V_z)$.
$$
\xymatrix{
&K_{p-2,q+1}(I_{S^q(X_z)},V_z) \ar[d] \\
K_{p-1,q+1}(I_{S^q(X)},V) \ar[r]_{\text{pr}_z} \ar@{-->}[ru] & K_{p-2,q+1}(I_{S^q(X)},V_z)
}
$$
This result need not be true anymore if we do not impose the condition $z\in X$.
\end{Rmk}

\begin{Cor}\label{Koszulrank}
Suppose that $K_{p,q}(S^q(X))\neq 0$. Then, for any nonzero Koszul class $\eta$ lying in $K_{p-1,q+1}(I_{S^q(X)},V)$, the induced map $\eta:\bigwedge^{p-1}V^\ast\to(I_{S^q(X)})_{q+1}$ has rank at least $\binom{p+q}{q+1}$. In particular, we have
$$
\dim(I_{S^q(X)})_{q+1}\geq\binom{p+q}{q+1}.
$$
\end{Cor}

\begin{proof}
 Let $\eta$ be a nonzero Koszul class in $K_{p-1,q+1}(I_{S^q(X)},V)$. Choose $r+1$ general points $z_0,\ldots,z_r\in V^\ast$ in $X$ with their dual basis $x_0,\ldots,x_r\in V$.
The collection of
$$
f_{\{0,\ldots,p+q-1\}\setminus\beta},\quad|\beta|=q+1\text{ in }\{0,\ldots,p+q-1\}
$$
gives the desired result by the proposition above. Here, $\eta\neq0$ makes sure that $p+q-1\leq r$.
\end{proof}

The main result of this section, Theorem \ref{kp1}, is now proved as follows:

\begin{proof}[\bf Proof of Theorem \ref{kp1}]
Suppose that $K_{p,q}(S^q(X))\neq 0$ with $p\geq e$. Let $X_\Lambda$ be a general $(e-1)$-inner projection of $X$. By Proposition \ref{monomial}, we obtain $K_{p-e+1,q}(S^q(X_\Lambda))\neq 0$. But since $S^q(X_\Lambda)$ is a hypersurface of degree larger than $q$, we have $p=e$ and $\deg S^q(X_\Lambda)=q+1$. By Theorem \ref{gBcwf}, $S^q(X)$ is of minimal degree.

Now, for the reverse implication, assume that $S^q(X)$ is a $q$-secant variety of minimal degree. We will show that $\text{reg}S^q(X)=q+1$ instead, which is enough since the projective dimension of $S_{S^q(X)}$ is at least $e$. Keep in mind that the case $e=1$ is trivial and the case $q=1$ is well known, and that $\text{codim}S^q(X_z)=e-1$ and $\text{codim}S^{q-1}(X_{\mathbb T_zX})=e$ for a general point $z\in X$. We are done by double induction on $e,q$ together with Corollary \ref{degbound}, Proposition \ref{paring} and Theorem \ref{vanishing}.
\end{proof}

\section{Higher secant varieties of minimal degree}

\subsection{Syzygetic characterization of minimal degree $q$-secant varieties}

For convenience, one can say that the Betti table of $X\subseteq\mathbb P^r$ is {\it $q$-pure Cohen-Macaulay} if it is of the following form.
$$
\begin{array}{c|cccc}
       & 0 & 1     & \cdots & e     \\ \hline
0 & 1 &            &                &            \\
q &        & B^e_{1,q} & \cdots & B^e_{e,q}
\end{array}
$$
Here, the remaining entries are all zero, and as in Corollary \ref{syzbound},
$$
B^e_{p,q}=\binom{p+q-1}{q}\binom{e+q}{p+q}.
$$

Due to Proposition \ref{basic} (or Proposition \ref{alternating}) and Theorem \ref{vanishing}, we have the following lemma:
\begin{Lem}\label{qpure}
Let $e=\textup{codim}S^q(X)$. Set $e\geq 2$ and $q\geq 2$. Assume that $S^q(X)$ has simple tangent cone at a general point $z$ of $X$. If
\begin{enumerate}[\normalfont (1)]
\item $S^q(X_z)$ has $q$-pure Cohen-Macaulay Betti table; and

\item $S^{q-1}(X_{\mathbb T_zX})$ has $(q-1)$-pure Cohen-Macaulay Betti table,
\end{enumerate}
then $S^q(X)$ has $q$-pure Cohen-Macaulay Betti table.
\end{Lem}

Now, we proceed to prove one of our main theorems.

\begin{proof}[\bf Proof of Theorem \ref{syzM}]
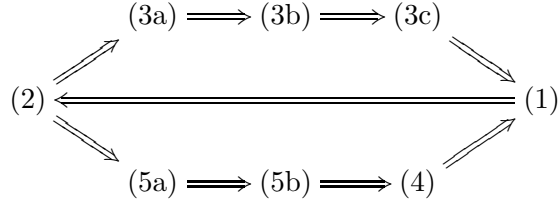
$$
\xymatrix@R=5mm{
 & \text{(3a)} \ar@{=>}[r] & \text{(3b)} \ar@{=>}[r] & \text{(3c)} \ar@{=>}[rd] & \\
\text{(2)} \ar@{=>}[ru] \ar@{=>}[rd] & & & & \text{(1)} \ar@{=>}[llll] \\
 & \text{(5a)} \ar@{=>}[r] & \text{(5b)} \ar@{=>}[r] & \text{(4)} \ar@{=>}[ru] &
}
$$

The directions $\text{(2)}\Rightarrow\text{(3a)}\Rightarrow\text{(3b)}\Rightarrow\text{(3c)}$ and $\text{(2)}\Rightarrow\text{(5a)}\Rightarrow\text{(5b)}\Rightarrow\text{(4)}$ are trivial. The direction $\text{(4)}\Rightarrow\text{(1)}$ is by the generalized $K_{p,1}$ theorem (Theorem \ref{kp1}). It suffices to show $\text{(3c)}\Rightarrow\text{(1)}$ and $\text{(1)}\Rightarrow\text{(2)}$.

In order to prove $\text{(3c)}\Rightarrow\text{(1)}$, suppose that $\text{(3c)}$ holds. Then, the $q$-secant variety to a general $(e-1)$-inner projection of $X$ has minimal degree by the equality case in Corollary \ref{syzbound}. Thus, by Theorem \ref{gBcwf}, $S^q(X)$ is of minimal degree.

For a proof of $\text{(1)}\Rightarrow\text{(2)}$, mimic the second part in the proof of Theorem \ref{kp1}. One shall use Lemma \ref{qpure}.
\end{proof}

\section{Del Pezzo higher secant varieties}
Already mentioned in the introduction, {\it del Pezzo higher secant varieties} are subextremal objects in the family of higher secant varieties. This section provides their definition and several properties that look similar to those of minimal degree higher secant varieties.

\begin{Prop}\label{degboundII}
Let $e=\textup{codim}S^q(X)$. If $S^q(X)$ is not a $q$-secant variety of minimal degree, then
$$
\deg S^q(X)\geq\binom{e+q}{q}+\binom{e+q-1}{q-1}=:D'_{e,q}.
$$
And when $e\geq 2$, the equality holds if and only if for a general point $z$ in $X$,
$$
\deg S^q(X_z)=D'_{e-1,q}\quad\text{and}\quad\deg S^{q-1}(X_{\mathbb T_zX})=D'_{e,q-1}
$$
together with the map $S^q(X)\dashrightarrow S^q(X_z)$ birational and $\deg\mathbb TC_zS^q(X)=\deg S^{q-1}(X_{\mathbb T_zX})$.
\end{Prop}

\begin{proof}
First of all, note that by Theorem \ref{gBcwf}, if $S^q(X)$ does not have minimal degree, then neither does $S^q(X_z)$ (resp., $S^{q-1}(X_{\mathbb T_zX})$) whenever $e\geq 2$ (resp., $q\geq 2$).

Let us see the case $e=1$. We will use induction on $q\geq 1$. Let $s=\deg\{\pi_z:S^q(X)\dashrightarrow S^q(X_z)\}$ for a general point $z\in X$. If $s=1$, then one can easily find that $\deg S^q(X)=q+1$, a contradiction. So, we have $\deg S^q(X)\geq 2+\deg S^{q-1}(X_{\mathbb T_zX})$ by the degree formula \eqref{degformula}. And by Theorem \ref{gBcwf}, $S^{q-1}(X_{\mathbb T_zX})$ do not have minimal degree so that the induction process works. On the other hand, the case $q=1$ is trivial. For the remaining cases, use double induction on $e$ and $q$, together with the Theorem \ref{gBcwf}. Refer to the identity $D'_{e-1,q}+D'_{e,q-1}=D'_{e,q}$.
\end{proof}

\begin{Rmk}
Let $e=\textup{codim}S^q(X)$. Suppose that $\deg S^q(X)=\binom{e+q}{q}+\binom{e+q-1}{q-1}$, and that the map $\widetilde{S}^q(X)\to S^q(X)$ is birational. Then, a general $e$-inner projection of $X$ is a {\it variety of two apparent $q$-secant $(q-2)$-planes}.
\end{Rmk}

\begin{Rmk}\label{positive}
In positive characteristic, Proposition \ref{degboundII} can not apply. For example, with $k$ equal to the algebraic closure of ${\mathbb F}_2$, consider a rational map
$\mathbb P_k^4\dashrightarrow\mathbb P_k^8$ given by the linear system
$$
\langle y_0y_2^2,y_0y_3^2,y_0y_4^2,y_1y_2^2,y_1y_3^2,y_1y_4^2,y_2^3,y_2^2y_3,y_2^2y_4\rangle\subset|\mathcal{O}_{\mathbb P_k^4}(3)|.
$$
Its image $X\subset\mathbb P_k^8$ is a nondegenerate irreducible projective variety whose homogeneous ideal is generated by the $2$-minors of the matrix
$$
M=
\begin{pmatrix}
x_0 & x_1 & x_2 \\
x_3 & x_4 & x_5 \\
x_6^2 & x_7^2 & x_8^2
\end{pmatrix}.
$$
Then, thanks to the identity $M(z+z')=M(z)+M(z')$, $z,z'\in V^\ast$ (Frobenius in characteristic two), the $2$-secant variety $S^2(X)$ is defined by the determinant of $M$. Hence, it has degree $4$, whereas the inequality above says that $\deg S^2(X)\geq 5$.

In characteristic zero, the projective variety $X\subset\mathbb P^8$ constructed in the same way satisfies $S^2(X)=\mathbb P^8$.
\end{Rmk}

The {\it sectional genus} $\pi(X)$ of $X$ is the arithmetic genus of a general curve section of $X$. So, the Hilbert polynomial of an $n$-dimensional projective variety $X$ is given by
$$
\deg X\cdot\binom{z+n-1}{n}+(1-\pi(X))\binom{z+n-2}{n-1}+\text{(lower terms)}.
$$

By Theorem \ref{secgenus}, we have the following lemma:

\begin{Lem}
 Suppose that the equality in Proposition \ref{degboundII} holds, and that $\textup{codim}S^q(X)\geq 2$. Then, we obtain
$$
\pi(S^q(X))\leq \pi(S^q(X_z))+\pi(S^{q-1}(X_{\mathbb T_zX}))+\deg S^{q-1}(X_{\mathbb T_zX})-1.
$$
\end{Lem}

\begin{proof}
Following the notation in Theorem \ref{secgenus}, one may see that the projective schemes $Z_i$ contain $\mathbb PTC_zS^q(X)$, and have the same dimension as $\mathbb PTC_zS^q(X)$. So, we have $\deg Z_i\geq\deg\mathbb PTC_zS^q(X)$. By these inequalities, we are done.
\end{proof}

By using the lemma above, double induction on $e,q$ deduces the following:

\begin{Prop}\label{secgenbound}
Let $e=\textup{codim}S^q(X)$. Suppose that $\deg S^q(X)=\binom{e+q}{q}+\binom{e+q-1}{q-1}$. Then, we obtain
$$
\pi(S^q(X))\leq(q-1)\left(\binom{e+q}{q}+\binom{e+q-1}{q-1}\right)+1.
$$
And if the equality holds for the case $e\geq 2$, then for a general point $z\in X$, $\pi(S^q(X_z))$ and $\pi(S^q(X_{\mathbb T_zX}))$ have their own maximal values.
\end{Prop}

\begin{Rmk}\label{secgenM}
Similarly, if $S^q(X)\subseteq\mathbb P^r$ is a $q$-secant variety of minimal degree, then
$$
\pi(S^q(X))=(q-1)\binom{e+q}{q}-\binom{e+q}{q-1}+1.
$$
Here, we used the equality instead of the inequality because we have $s=s'=1$ in Theorem \ref{secgenus}.
\end{Rmk}

\begin{Def}\label{dP}
A $q$-secant variety of {\it almost minimal degree} is defined to be such that for its degree, the equality in Proposition \ref{degboundII} holds. And a $q$-secant variety $S^q(X)$ of almost minimal degree is called a {\it del Pezzo $q$-secant variety} if $\pi(S^q(X))$ attains the maximum in Proposition \ref{secgenbound}.
\end{Def}

\begin{Rmk}\label{dPfunda}
Here are some fundamental remarks:
\begin{enumerate}[\normalfont (1)]
\item Let $X\subset\mathbb P^r$ be of almost minimal degree. Then, as aforementioned in the introduction, we have
$$
\pi(X)=
\begin{cases}
1 & \text{if $X$ is a del Pezzo variety} \\
0 & \text{otherwise.}
\end{cases}
$$

\item Suppose that $S^q(X)$ is del Pezzo of codimension $e$, and let $z\in X$ be a general point. If $e\geq 2$ and $q\geq 1$, then $S^q(X_z)$ is del Pezzo. If $e\geq 1$ and $q\geq 2$, then $S^{q-1}(X_{\mathbb T_zX})$ is del Pezzo.
\end{enumerate}
\end{Rmk}

\noindent{\bf Examples of del Pezzo higher secant varieties.}
\begin{enumerate}[\normalfont (1)]
\item Let $E\subset\mathbb P^r$ be an elliptic normal curve. Then, one may find that the $S^q(E)$ with $e:=\text{codim}S^q(E)>0$ are del Pezzo with Betti table
$$
\begin{tabularx}{0.5\linewidth}{ c| *{5}{>{\centering\arraybackslash}X} }
         & $0$ & $1$     & $\cdots$ & $e-1$ & $e$ \\ \hline
$0$   & $1$ &            &               &          & \\
$q$   &        & $\ast$ & $\cdots$ & $\ast$ & \\
$2q$ &        &            &             &           & $1$
\end{tabularx}
$$
Refer to \cite{bothmer2004geometric} and \cite{fisher2006higher}. In fact, Proposition \ref{gBcIIcurve} and Example \ref{dPcurve} say that a nondegenerate irreducible projective curve $C\subset\mathbb P^r$ with $e:=\text{codim}S^q(C)\geq 2$ is del Pezzo if and only if $S^q(C)$ is del Pezzo.

\item Let $Y\subset\mathbb P^{11}$ be the complete embedding of $S(1,2)\subset\mathbb P^4$ by $|\mathcal O_{S(1,2)}(2)|$. Then, using Macaulay2(\cite{M2}), we find that $S^2(Y)$ is del Pezzo of codimension $6$. Moreover, by Theorem \ref{gBcIIwf}, one will be able to see that if $X=\nu_4(\mathbb P^2)\subset\mathbb P^{14}$ is the fourth Veronese surface, then $S^3(X)$ is also del Pezzo of codimension $6$ because a general tangential projection of $X$ is projectively equivalent to $Y$. (Note that $S^4(X)$ has minimal degree and codimension $3$.)

\item For a higher dimensional example, consider the Pl{\"u}cker embedding $X=\mathbb G(\mathbb P^1,\mathbb P^{2q+2})\subset\mathbb P^r$, $r=\binom{2q+3}{2}-1$, of the Grassmannian. Then, it is defined by the ideal $\text{Pf}_4(M)$ for the $(2q+3)\times(2q+3)$ generic skew-symmetric matrix $M$, and thus the $q$-secant variety $S^q(X)$ by the ideal $\text{Pf}_{2q+2}(M)$. Apparently, $S^q(X)$ has the same Betti table as $S^q(E)\subset\mathbb P^{2q+2}$ in the paragraph above (see \cite{buchsbaum1977algebra}), and is a del Pezzo $q$-secant variety. Refer to the fact that a general tangential projection of $\mathbb G(\mathbb P^1,\mathbb P^{2q+2})$ is $\mathbb G(\mathbb P^1,\mathbb P^{2q})\subset\mathbb P^{r'}$, $r'=\binom{2q+1}{2}-1$. This is because by symmetry, a general point $z\in V^\ast$ in $\mathbb G(\mathbb P^1,\mathbb P^{2q+2})$ can be assumed to satisfy $M(z)_{0,1}=1$, $M(z)_{1,0}=-1$ and $M(z)_{i,j}=0$ otherwise, and therefore, the submatrix of $M$ obtained by deleting first two rows and columns defines $\mathbb G(\mathbb P^1,\mathbb P^{2q})$.
\end{enumerate}

\begin{Rmk}\label{ci}
If $S^q(X)$ is a del Pezzo $q$-secant variety of codimension two, then it is a complete intersection of two $(q+1)$-forms. Indeed, the case $q=1$ is classical, and if $q=2$, then applying \cite[p.\phantom{.}32]{gruson1978genre} and \cite[p.\phantom{.}251]{roggero2003laudal} to its general curve section $C:=S^2(X)\cap\Lambda$, we find that $C$ is a complete intersection of two cubics, hence so is $S^2(X)$. And for the other cases $q\geq 3$, use Lemma \ref{prolongation} to get two linearly independent $(q+1)$-forms; they generate the homogeneous ideal of $S^q(X)$ due to the invariants of $S^q(X)$.
\end{Rmk}

Note that by \cite{han2015sharp}, if $X$ is not a variety of minimal degree with codimension $e$, then we have $\beta_{p,1}(X)\leq p\binom{e+1}{p+1}-\binom{e}{p-1}$, $1\leq p\leq e$. By double induction on $e,q$ based on Theorem \ref{gBcwf}, the following is straightforward:

\begin{Prop}\label{syzboundII}
Let $e=\textup{codim}S^q(X)$. Suppose that $S^q(X)$ is not a $q$-secant variety of minimal degree. Then,
$$
\beta_{p,q}(S^q(X))\leq B'^e_{p,q}:=\binom{p+q-1}{q}\binom{e+q}{p+q}-\binom{e+q-p-1}{q-1}\binom{e+q-1}{e+q-p}
$$
for every $1\leq p\leq e$. And if the equality holds for some $1\leq p\leq e$, then for a general point $z$ in $X$, we have $
\beta_{p,q}(S^q(X_z))=B'^{e-1}_{p,q}$, $\beta_{p-1,q}(S^q(X_z))=B'^{e-1}_{p-1,q}$ and $\beta_{p,q-1}(S^{q-1}(X_{\mathbb T_zX}))=B'^e_{p,q-1}$.
\end{Prop}

\subsection{An analogue of the generalized Bronowski's conjecture}

One says $X$ to be a {\it $\textup{dP}^q$-variety} if $S^q(X)$ is a del Pezzo $q$-secant variety, and to be a {\it $\textup{dPA}^q$-variety} if it is a $\text{dP}^q$-variety and the map $\widetilde{S}^q(X)\to S^q(X)$ is birational. In view of the generalized Bronowski's conjecture, one might expect the following:

\begin{Conj}[The {\it generalized Bronowski's conjecture II}]\label{gBcII}
If a general $(q-1)$-tangential projection of $X$ is birational onto a del Pezzo variety of codimension at least $2$, then $X$ is a $\textup{dPA}^q$-variety.
\end{Conj}

Note that the condition on the codimension can not be dropped because for a smooth rational curve $X=C\subset\mathbb P^4$ of degree $5$, its general tangential projection is birational onto a plane cubic, but $\deg S^2(C)=6\neq 5$. Also, see the following:

\begin{Prop}\label{gBcIIconverse}
The converse of the generalized Bronowski's conjecture II holds, that is, if $X$ is a $\textup{dPA}^q$-variety with $\textup{codim}S^q(X)\geq 2$, and if $\pi_\Lambda:X\dashrightarrow X_\Lambda$ is a general $(q-1)$-tangential projection, then $\pi_\Lambda$ is birational and $X_\Lambda$ is a del Pezzo variety.
\end{Prop}

\begin{proof}
By Remark \ref{dPfunda}, $X_\Lambda$ is a del Pezzo variety. It suffices to show that $\pi_\Lambda:X\dashrightarrow X_\Lambda$ is birational.

Let $e$ be the codimension of $S^q(X)$, and $\Lambda'$ the span of $e$ general points in $X$. It is easy to see that the map $\widetilde{S}^q(X_{\Lambda'})\to S^q(X_{\Lambda'})$ has degree two since $S^q(X)$ has almost minimal degree with the map $\widetilde{S}^q(X)\to S^q(X)$ birational. By \cite[Theorem 2.7]{ciliberto2006varieties}, we have
$$
\deg\{\pi_\Lambda:X_{\Lambda'}\to X_{\langle\Lambda,\Lambda'\rangle}\}\leq 2.
$$
Further, note that $X_\Lambda$ has almost minimal degree, hence $\pi_{\Lambda'}:X_\Lambda\to X_{\langle\Lambda,\Lambda'\rangle}$ has degree two. Consider the following commutative diagram:
$$
\xymatrix{
X_{\langle\Lambda,\Lambda'\rangle}&X_\Lambda \ar@{-->}[l]_--{2:1}\\
X_{\Lambda'} \ar@{-->}[u] &X \ar@{-->}[u] \ar@{-->}[l]^{1:1}
}
$$
We find that $\pi_\Lambda:X_{\Lambda'}\to X_{\langle\Lambda,\Lambda'\rangle}$ has degree two and $\pi_\Lambda:X\dashrightarrow X_\Lambda$ is birational.
\end{proof}

Referring to Proposition \ref{paring} and its proof, together with Remark \ref{ci}, we have:

\begin{Prop}\label{paringII}
Set $e:=\textup{codim}S^q(X)\geq 2$. If $S^q(X)$ is a del Pezzo $q$-secant variety, then it satisfies property $ST$, and if in addition, $e\geq 3$, then it satisfies property $IC$.
\end{Prop}

The following is a proof of a weak form of the generalized Bronowski's conjecture II:

\begin{proof}[\bf Proof of Theorem \ref{gBcIIwf}]
(1) Set $e\geq3$. We will use induction on $q\geq 1$. If a general inner projection of $X$ is del Pezzo, then by Proposition \ref{paringII} and Proposition \ref{paring}(1), $X$ has simple tangent cone at its general point and thus is del Pezzo; we have dealt with the case $q=1$. Now, assume that $q\geq 2$. Let $w$ be another general point of $X$. Then, $S^{q-1}(X_{\langle w,\mathbb T_zX\rangle})$ is a del Pezzo $(q-1)$-secant variety. By the induction hypothesis, $S^{q-1}(X_{\mathbb T_zX})$ is a del Pezzo $(q-1)$-secant variety. Proposition \ref{paring} applies so that $S^q(X)$ has simple tangent cone at a general point of $X$. Therefore, $S^q(X)$ is a del Pezzo $q$-secant variety.

(2) Similarly, use induction on $e\geq 2$. Refer to the proof of (2) in Theorem \ref{gBcwf}.
\end{proof}

\begin{Prop}\label{gBcIIcurve}
Let $X=C\subset\mathbb P^r$, $r\geq 5$, be a curve. Suppose that its general tangential projection is a del Pezzo curve. Then, $S^2(C)$ is a del Pezzo secant variety. Furthermore, $C$ is a del Pezzo curve. Especially, the generalized Bronowski's conjecture II is true for curves.
\end{Prop}

The proof will be given later after establishing lemmas.

\begin{Ex}\label{dPcurve}
If $C\subset\mathbb P^r$ is a del Pezzo curve, then $S^q(C)$ is a del Pezzo $q$-secant variety whenever $S^q(C)\neq\mathbb P^r$. Indeed, suppose that $C$ is a del Pezzo curve with $S^q(C)\neq\mathbb P^r$, and let us prove that $S^q(C)$ is a del Pezzo $q$-scant variety. If $C$ is an elliptic normal curve, then we are done by \cite{fisher2006higher}. Assume that $C$ is singular so that $C$ is the projection of a rational normal curve $D\subset\mathbb P^{r+1}$ from a point $z\in S^2(D)\setminus D$. By Theorem \ref{gBcwf}, taking inner projections several times, it suffices to consider the case $r=(2q-1)+2$ only. We will mimic \cite{lee2013non}. Note that the homogeneous ideal of $S^q(D)$ is generated by $(q+1)$-minors of the matrix
$$
M=
\begin{pmatrix}
x_0 & x_1 & \cdots & x_q & x_{q+1} & x_{q+2} \\
x_1 & x_2 & \cdots & x_{q+1}& x_{q+2} & x_{q+3} \\
\vdots & \vdots & \ddots & \vdots & \vdots & \vdots \\
x_q & x_{q+1} & \cdots & x_{2q} & x_{2q+1} & x_{2q+2}
\end{pmatrix}.
$$
By symmetry, we need to consider only two cases: $z=(1:0:\cdots:0:1)$ or $z=(0:1:0:\cdots:0)$. Let $\Delta^{i,j}$ be the minor of $M$ obtained by deleting $i$-th and $j$-th columns. If $z=(1:0:\cdots:0:1)$, then $\Delta^{0,q+2}$ and $\Delta^{q+1,q+2}-\Delta^{0,1}$ are linearly independent two cubic equations of $S^2(C)$. And if $z=(0:1:0:\cdots:0)$, then $\Delta^{0,1}$ and $\Delta^{1,2}+\Delta^{0,3}$ are. Thus, $S^q(C)$ is del Pezzo.
\end{Ex}

\begin{Lem}\label{tanalmost}
Let $X=C\subset\mathbb P^r$, $r\geq 5$, be a curve. If its general tangential projection has almost minimal degree, then so does itself.
\end{Lem}

\begin{proof}
Let $z$ be a general point of $C$. Then, either $C_{\mathbb T_zC}$ is smooth, or its homogeneous ideal is generated by quadrics. Thus, in any event, $C_{\mathbb T_zC}$ is not tangentially degenerate. (Refer to \cite{kaji1986tangentially} for smooth curves.) It is easy to see that $C$ is also not tangentially degenerate. We have $\deg C=\deg C_{\mathbb T_zC}+2$.
\end{proof}

\begin{Lem}\label{sing}
Suppose that the homogeneous ideal of $S^{q+1}(X)$ is generated by $(q+2)$-forms, and that $\textup{Sing}S^{q+1}(X)=S^q(X)$. Let $N=\dim S^{q+1}(X)$. Then, we have
$$
\bigcap_{i=0}^{N+1}J(\mathbb T_{z_i}X,S^q(X))=S^q(X).
$$
as sets for $N+2$ general points $z_0,\ldots,z_{N+1}\in X$, where $J$ means the \textup{join}.
\end{Lem}

\begin{proof}
Take $N+2$ general points $z_0,\ldots,z_{N+1}$ in $X$, and choose homogeneous coordinates $x_0,\ldots,x_r$ such that $z_i\in V(x_j)$ if and only if $i\neq j$. Let $p$ be any point of the intersection in question. Then, the first $N+2$ columns of the Jacobian matrix vanish at $p$. Therefore, $p\in\text{Sing}S^{q+1}(X)=S^q(X)$. The reverse containment is trivial.
\end{proof}

\begin{proof}[Proof of Proposition \ref{gBcIIcurve}]
By Lemma \ref{tanalmost}, $C$ is of almost minimal degree. To the contrary, assume that $C$ is an isomorphic projection of a rational normal curve $D\subset\mathbb P^{r+1}$ from a point $b\in\mathbb P^{r+1}\setminus S^2(D)$. Let $w$ be a general point of $D$ so that $z:=\pi_b(w)$ is also a general point of $C$. Let $\overline{b}=\pi_{\mathbb T_wD}(b)$. The tangential projection $D_{\mathbb T_wD}$ is obviously a rational normal curve in $\mathbb P^{r-1}$.
$$
\xymatrix{
C_{\mathbb T_zC} & D_{\mathbb T_wD} \ar@{-->}[l]_{\pi_{\overline{b}}} \\
C \ar@{-->}[u]^{\pi_{\mathbb T_zC}} & D \ar[l]^{\pi_b} \ar@{-->}[u]_{\pi_{\mathbb T_wD}}
}
$$
We claim that $\overline{b}\not\in S^2(D_{\mathbb T_wD})$. Since $w$ was general in $D$, by Lemma \ref{sing} together with \cite[p.\phantom{.}440]{bertram1992moduli}, $b\not\in J(\mathbb T_wD,S^2(D))$, or  in other words, $\overline{b}\not\in S^2(D_{\mathbb T_wD})$. Thus, $\pi_{\overline{b}}$ is not only a regular map but also an isomorphic projection, a contradiction. Thus, $C$ is a del Pezzo curve. Now, in order to show that $S^2(C)$ is a del Pezzo secant variety, refer to Example \ref{dPcurve}.
\end{proof}

\begin{Rmk}
Let $q\geq 2$, and assume that $\text{codim}S^q(C)\geq 2$. If $S^q(C)$ has almost minimal degree, then it is in fact a del Pezzo $q$-secant variety. Indeed, let $C'$ be a general $(q-2)$-tangential projection of $C$ so that $S^2(C')$ has almost minimal degree. Since a general tangential projection of $C'$ is of almost minimal degree, so is $C'$. If $C'$ were an isomorphic projection of a rational normal curve, then $\deg S^2(C')=\binom{(e+1)+2}{2}\neq\binom{e+2}{2}+\binom{e+1}{1}$ for $e:=\text{codim}S^2(C')$. Thus, $C'$ is a del Pezzo variety, and $S^2(C')$ is a del Pezzo $2$-secant variety by Example \ref{dPcurve}. By Theorem \ref{gBcIIwf}, $S^q(C)$ is a del Pezzo $q$-secant variety.
\end{Rmk}

\subsection{Syzygetic characterization of del Pezzo higher secant varieties}

We need the technical lemma below:

\begin{Lem}\label{technical}
Let $I\subseteq S$ be a homogeneous ideal with $I_q=0$. Suppose that $S/I$ satisfies property $N_{q+1,p}$, $p\geq 1$. Let $m=\min\{j:\beta_{p+1,j}(S/I)\neq 0\}$ with $k=p+m+1$ and $m'=m-q+1$. Consider a double complex
$$
C^{i,j}=\bigwedge^{k-i-j}V\otimes I_i\otimes S_j
$$
with differentials $\delta_I,\delta_{II}$ induced by Koszul type differentials. Then, there is an isomorphism
$$
K_{p+1,m}(S/I,V)\cong(\ker\delta_I^{q+1,m'})\cap(\ker\delta_{II}^{q+1,m'})\subseteq C^{q+1,m'}.
$$
And the intersection above maps to zero under the multiplication map 
$$
C^{q+1,m'}=\bigwedge^{p-1}V\otimes I_{q+1}\otimes S_{m'}\to\bigwedge^{p-1}V\otimes I_{m+2}.
$$
\end{Lem}

\begin{proof}
Since $S/I$ satisfies property $N_{q+1,p}$, chopping the minimal free resolution of $S/I$ into short exact sequences, we have graded $S$-modules $I=Q_1,Q_2,\ldots,Q_{p+1}$ with
$$
\xymatrix{
0 \ar[r] & Q_{i+1} \ar[r] & (Q_i)_{q+i}\otimes S(-q-i) \ar[r] & Q_i \ar[r] & 0
}
$$
for each $1\leq i\leq p$. Then, it is straightforward to show that $K_{p+1,m}(S/I)\cong K_{0,k}(Q_{p+1})\cong (Q_{p+1})_k$. Note that
$$
(Q_{p+1})_k=\ker\{(Q_p)_{p+q}\otimes S_{m-q+1}\to (Q_p)_k\},
$$
$$
(Q_p)_{p+q}=\ker\{(Q_{p-1})_{p+q-1}\otimes V\to (Q_{p-1})_{p+q}\},
$$
and so on. That is, $(Q_{p+1})_k\subseteq I_{q+1}\otimes T^{p-1}V\otimes S_{m'}$ consists of elements mapping to zero under the maps induced by multiplication maps of adjacent pairs. Considering multiplication maps of adjacent pairs in $T^{p-1}V$, one can see that the inclusion $(Q_{p+1})_k\subseteq I_{q+1}\otimes T^{p-1}V\otimes S_{m'}$ factors through $I_{q+1}\otimes \bigwedge^{p-1}V\otimes S_{m'}=\bigwedge^{p-1}V\otimes I_{q+1}\otimes S_{m'}=C^{q+1,m'}$. And in consideration of the multiplication maps $I_{q+1}\otimes V\to I_{q+2}$ and $V\otimes S_{m'}\to S_{m'+1}$ of adjacent pairs, the first assertion holds.

 Now, let us prove the second assertion. Take any element $\theta\in(\ker\delta_I)\cap(\ker\delta_{II})\subseteq C^{q+1,m'}$. Consider maps $D:C^{i,j}\to C^{i,j-1}$ induced by $\omega\otimes f\mapsto\sum_{l=0}^r x_l\wedge\omega\otimes\partial f/\partial x_l$, and the element $(\delta_ID)^{m'}\theta\in C^{m+2,0}=\bigwedge^{p-1}V\otimes I_{m+2}$. Notice that
$$
D\delta_{II}+\delta_{II}D=(k-i)\cdot\text{id}
$$
on each $C^{i,j}$. Then, by (finite) induction, we see that $\delta_{II}(\delta_ID)^s\theta=0$ for all $0\leq s\leq m'$. But the map $\delta_{II}^{m+2,0}:\bigwedge^{p-1}V\otimes I_{m+2}\to\bigwedge^{p-2}V\otimes I_{m+2}\otimes V$ is injective, hence $(\delta_ID)^{m'}\theta=0$. On the other hand, $(\delta_ID)^{m'}\theta$ is nothing but the image of $\theta$ under the multiplication map $C^{q+1,m'}=\bigwedge^{p-1}V\otimes I_{q+1}\otimes S_{m'}\to\bigwedge^{p-1}V\otimes I_{m+2}$, up to nonzero factor. We are done.
\end{proof}

For convenience, one can say that the Betti table of $X\subseteq\mathbb P^r$ is {\it $q$-pure Gorenstein} if it is of the following form:
$$
\begin{tabularx}{0.6\linewidth}{ c| *{5}{>{\centering\arraybackslash}X} }
         & $0$ & $1$     & $\cdots$ & $e-1$ & $e$ \\ \hline
$0$   & $1$ &            &                &           & \\
$q$   &        & $B'^e_{1,q}$ & $\cdots$ & $B'^e_{e-1,q}$ & \\
$2q$ &        &            &                &           & $1$
\end{tabularx}
$$
Here, $e:=\text{codim}X$, the remaining entries are all zero, and the value of $B'^e_{p,q}$ is given in Proposition \ref{syzboundII}, namely
$$
B'^e_{p,q}=\binom{p+q-1}{q}\binom{e+q}{p+q}-\binom{e+q-p-1}{q-1}\binom{e+q-1}{e+q-p}.
$$

\begin{Lem}\label{pureGorenstein}
Let $e=\textup{codim}S^q(X)$. Set $e\geq 2$ and $q\geq 2$. Assume that $S^q(X)$ has simple tangent cone at a general point $z$ of $X$. If
\begin{enumerate}[\normalfont (1)]
\item $S^q(X_z)$ has $q$-pure Gorenstein Betti table; and

\item $S^{q-1}(X_{\mathbb T_zX})$ has $(q-1)$-pure Gorenstein Betti table,
\end{enumerate}
then $S^q(X)$ has $q$-pure Gorenstein Betti table.
\end{Lem}

\begin{proof}
By Theorem \ref{vanishing} and Proposition \ref{alternating}, all Betti numbers of $S^q(X)\subseteq\mathbb P^r$ have the required values except $\beta_{e-1,2q}(S^q(X))$ and $\beta_{e,2q-1}(S^q(X))$ whose values are the same. By considering the ``partial'' Euler characteristic of $\mathbb F\otimes_SK(S)$, where $\mathbb F$ is the minimal free resolution of $S_{S^q(X)}$, we find that the Betti number $\beta_{e-1,2q}(X)$ is zero. We are done.
\end{proof}

Let us prove the syzygetic characterization of del Pezzo higher secant varieties.

\begin{proof}[\bf Proof of Theorem \ref{syzdP}]
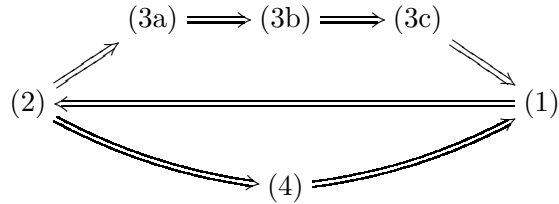
$$
\xymatrix@R=5mm{
& \text{(3a)} \ar@{=>}[r] & \text{(3b)} \ar@{=>}[r] & \text{(3c)} \ar@{=>}[rd] & \\
\text{(2)} \ar@{=>}[ur] \ar@/_/@{=>}[drr] & & & & \text{(1)} \ar@{=>}[llll] \\
 & & \text{(4)} \ar@/_/@{=>}[urr] & &
}
$$

The directions $\text{(2)}\Rightarrow\text{(3a)}\Rightarrow\text{(3b)}\Rightarrow\text{(3c)}$ and $\text{(2)}\Rightarrow\text{(4)}$ are trivial. So, we only need to show $\text{(3c)}\Rightarrow\text{(1)}$,$\text{(4)}\Rightarrow\text{(1)}$ and $(1)\Rightarrow (2)$.

To see $\text{(3c)}\Rightarrow\text{(1)}$, by Theorem \ref{gBcIIwf}, it suffices to show that $\beta_{1,q}(S^q(X_\Lambda))=2$ for a general $(e-2)$-inner projection $X_\Lambda$ of $X$. Indeed, it is the case thanks to the equality part of Proposition \ref{syzboundII}.

For $\text{(4)}\Rightarrow\text{(1)}$, consider a general $(e-2)$-inner projection $X_\Lambda$. We will show that $S^q(X_\Lambda)$ lies in two or more distinct hypersurfaces of degree $q+1$, which implies $\text{(1)}$ by Theorem \ref{gBcIIwf}. Notice that Proposition \ref{monomial} ensures the existence of at least one hypersurface of degree $q+1$ in which $S^q(X_\Lambda)$ is. Suppose that a unique degree $q+1$ hypersurface contains $S^q(X_\Lambda)$. Let $\eta\in K_{e-2,q+1}(I_{S^q(X)},V)$ be a nonzero Koszul class. Take $r+1$ general points $z_i\in V^\ast$ in $X$ with their dual basis $x_i$. Then, by Proposition \ref{monomial},
$$
\eta=\sum_{|\alpha|=e-2}x_\alpha\otimes f_\alpha
$$
with $f_\alpha$ nonzero in $I_{S^q(X_{\Lambda_\alpha})}$, where for any $\alpha=\{i_1<\cdots<i_{|\alpha|}\}\subseteq\{0,\ldots,r\}$, $x_\alpha=x_{i_1}\wedge\cdots\wedge x_{i_{|\alpha|}}\in\bigwedge^{|\alpha|}V$ and $\Lambda_\alpha=\langle z_i:i\in\alpha\rangle$. Let $\eta'\in K_{e-2,q+1}(I_{S^q(X)},V)$ be another Koszul class so that $\eta'=\sum_\alpha x_\alpha\otimes f'_\alpha$ with $f'_\alpha\in I_{S^q(X_{\Lambda_\alpha})}$ by Proposition \ref{monomial}. By assumption, $f'_\alpha$ is a multiple of $f_\alpha$.

 Now, referring to Lemma \ref{technical}, we will show that $\ker\delta_I\cap\ker\delta_{II}=0$ in $C^{q+1,m'}$ with $p=e-1$, which deduces that the projective dimension of $S_X$ is less than $e$, a contradiction. Let $\theta\in\ker\delta_I\cap\ker\delta_{II}\subseteq C^{q+1,m'}$ be an arbitrary element. By the observation above, we can write
$$
\theta=\sum_{|\alpha|=e-2}x_\alpha\otimes f_\alpha\otimes g_\alpha
$$
for some $g_\alpha\in S_{m'}$. But by Lemma \ref{technical} again, we have $\sum_\alpha x_\alpha\otimes f_\alpha g_\alpha=0$, which means that $g_\alpha=0$ for every $|\alpha|=e-2$, hence $\theta=0$.

Now, let us see the direction $(1)\Rightarrow (2)$. If $q=1$, then it is well known. Assume that $q\geq 2$. When $X_\Lambda$ is a general $(e-2)$-inner projection of $X$, consider its $q$-secant variety $S^q(X_\Lambda)$ whose codimension is $2$. By Remark \ref{ci}, it is a complete intersection of two $(q+1)$-forms. We are done by using double induction on $e,q$ with Proposition \ref{paringII} and Lemma \ref{pureGorenstein}.
\end{proof}

\begin{Rmk}
By \cite[Theorem 2.1]{buchsbaum1977algebra}, Theorem \ref{syzdP} implies that a del Pezzo $q$-secant variety of codimension $3$ is ideal-theoretically defined by the ideal $\text{Pf}_{2q+2}(M)$ of all the $(2q+2)$-{\it Pfaffians} of $M$ for some $(2q+3)\times(2q+3)$ skew-symmetric matrix $M$ of linear forms.

On the other hand, if there is a $(2q+3)\times(2q+3)$ skew-symmetric linear matrix $M$ such that its $4$-Pfaffians vanish on $X$ and $\text{codim}V(\text{Pf}_{2q+2}(M))=3$, then the $q$-secant variety $S^q(X)$ is del Pezzo with $I_{S^q(X)}=\text{Pf}_{2q+2}(M)$ unless it has minimal degree. This is because we have $\text{rank}M(w)\leq 2q$ for every $w\in S^q(X)$ (since $\text{rank}M(z)\leq 2$ on $X$) and thus $\dim(I_{S^q(X)})_{q+1}\geq\dim\text{Pf}_{2q+2}(M)_{q+1}=2q+3$ (due to $\text{codim}V(\text{Pf}_{2q+2}(M))=3$).
\end{Rmk}

\subsection{Higher secant varieties of almost minimal degree}

Let us remind that by definition, $S^q(X)$ is a $q$-secant variety of almost minimal degree if $\deg S^q(X)=\binom{e+q}{q}+\binom{e+q-1}{q-1}$, where $e:=\text{codim}S^q(X)$.

\begin{Prop}\label{almost}
Suppose that $S^q(X)$ is a $q$-secant variety of minimal degree with codimension $e\geq 2$ and $q\geq 2$. Then, $S^q(X_z)$ is a $q$-secant variety of almost minimal degree with codimension $e-1$ when $z$ is a general point of $S^2(X)$. In this case, $\deg\{\pi_z:S^q(X)\dashrightarrow S^q(X_z)\}=1$, $K_\infty(z,I_{S^q(X)})=(I_{S^{q-2}(X_\Lambda)})$ with $s(z,I_{S^q(X)})=2$, and
$$
\begin{aligned}
\pi(S^q(X_z)) & =(q-1)\left(\binom{e+q-1}{q}+\binom{e+q-2}{q-1}\right)+1 \\
&\phantom{=}-\left(\binom{e+q-1}{q-1}+\binom{e+q-2}{q-2}-\deg Z_1\right),
\end{aligned}
$$
where $\Lambda=\mathbb T_zS^2(X)$ and $Z_1\subseteq\mathbb P^{r-1}$ is the projective scheme defined by $K_1(z,I_{S^q(X)})$.
\end{Prop}

\begin{proof}
We may assume that there are two general points $z_0,z_1\in X$ such that $\langle z_0,z_1\rangle$ has $z$ as its general point, hence $z_0=(1:0:0:\cdots)$, $z_1=(0:1:0:\cdots)$ and $z=(1:1:0:\cdots)$ in a homogeneous coordinates system $x_i$. By Terracini's lemma, $\langle\mathbb T_{z_0}X,\mathbb T_{z_1}X\rangle=\mathbb T_zS^2(X)=\Lambda$, and so $X_\Lambda$ is a general $2$-tangential projection which is an $\text{M}^{q-2}$-variety.  As $S^q(X)$ and $S^{q-1}(X_{\mathbb T_{z_0}X})$ have simple tangent cone,
$$
K_1(z_0,I_{S^q(X)})=(I_{S^{q-1}(X_{\mathbb T_{z_0}X})})\quad\text{and}\quad K_1(\overline{z_1},I_{S^{q-1}(X_{\mathbb T_{z_0}X})})=(I_{S^{q-2}(X_\Lambda)}),
$$
where $\overline{z_1}$ is the image of $z_1$ under the projection $\pi_{\mathbb T_{z_0}X}$, which means that the composition
$$
\xymatrix{
(I_{S^q(X)})_{q+1} \ar[r]^--{\partial/\partial x_0} & (I_{S^{q-1}(X_{\mathbb T_{z_0}X})})_q \ar[r]^--{\partial/\partial x_1} & (I_{S^{q-2}(X_\Lambda)})_{q-1}
}
$$
is surjective. Thus, every minimal generator $f'$ of $I_{S^{q-2}(X_\Lambda)}$ comes from an element $f=f' x_0x_1+g_0x_0+g_1x_1+h\in I_{S^q(X)}$ for some polynomials $g_0,g_1,h$ in $x_2,\ldots,x_r$. Consider a coordinate change: $y_i=x_i$ for $i\neq 1$ and $y_1=x_1-x_0$ for which $z=(1:0:\cdots:0)$. Then, in coordinates $y_i$, $f$ looks like $f'y_0(y_0+y_1)+g_0y_0+g_1(y_0+y_1)+h$ with $f',g_0,g_1,h$ in $y_2,\ldots,y_r$. On the other hand, there is a nonzero $(q+1)$-form $f$ in $I_{S^q(X_{z_1})}$ such that $\deg_{y_0}f=\deg_{x_0}f=1$. Thus,
$$
K_0(z,I_{S^q(X)})\subsetneq K_1(z,I_{S^q(X)})\subsetneq K_2(z,I_{S^q(X)})=K_\infty(z,I_{S^q(X)})
$$
with $K_\infty(z,I_{S^q(X)})=(I_{S^{q-2}(X_\Lambda)})$, that is, $s(z,I_{S^q(X)})=2$ and $\deg\pi_z=1$. Now, the degree formula \eqref{degformula} gives
$$
\begin{aligned}
\deg S^q(X_z) & =\frac{\deg S^q(X)-\deg\mathbb TC_zS^q(X)}{\deg\pi_z} \\
 & =\binom{e+q}{q}-\binom{e+q-2}{q-2} \\
 & =\binom{e+q-1}{q}+\binom{e+q-2}{q-1}.
\end{aligned}
$$
Further, by Theorem \ref{secgenus}, we obtain
$$
\pi(S^q(X))=\pi(S^q(X_z))+\pi(S^{q-2}(X_\Lambda))+2\deg S^{q-2}(X_\Lambda)-\deg Z_1-1.
$$
By using Remark \ref{secgenM} for $\pi(S^q(X))$ and $\pi(S^{q-2}(X_\Lambda))$, the computation is straightforward.
\end{proof}

\begin{Rmk}
Let us consider the $q$-secant variety $S^q(X_z)$ constructed in the proposition above. When $X$ is a rational normal curve, $S^q(X_z)$ is a del Pezzo $q$-secant variety by Example \ref{dPcurve}. Therefore, $\deg Z_1=\binom{e+q-1}{q-1}+\binom{e+q-2}{q-2}$ in this case. However, some computational evidence suggests that $S^q(X_z)$ is not del Pezzo in general.
\end{Rmk}

\begin{Ex}
In assistance of Macaulay2 (\cite{M2}), Proposition \ref{almost} produces an example of a non-del Pezzo $2$-secant variety of almost minimal degree. Let $X_z\subset\mathbb P^7$ be the projection of a rational normal surface $X\subset\mathbb P^8$ (for instance, $S(3,4)$) from a general point $z\in S^2(X)$. Then, $S^2(X_z)$ is a $2$-secant variety of almost minimal degree and codimension $2$, and the Betti table of $S^2(X_z)$ is given below:
$$
\begin{array}{c|cccc}
   & 0 & 1 & 2 & 3 \\ \hline
0 & 1 & -  & - & -  \\
1 & - & -  & - & - \\
2 & - & 1 & - & - \\
3 & - & - & - & - \\
4 & - & 6 & 9 & 3
\end{array}
$$
Moreover, the projective scheme $Z_1=V(K_1(z,I_{S^2(X)}))$ has degree one, and so the sectional genus of $S^2(X_z)$ is six by the formula in Proposition \ref{almost}, which can also be checked by Macaulay2.
\end{Ex}

\medskip

\noindent{\bf Some remarks and questions.} 
\begin{enumerate}[\normalfont (1)]
\item The matryoshka structure of higher secant varieties would include plenty of pictures in projective geometry other than what we have presented. For example, M. Green’s {\it $(2g+1+ p)$-theorem} is such a picture. The theorem says that for a linearly normal smooth curve $X=C\subseteq\mathbb P^r$ of genus $g$, if $C$ has degree $d\geq 2g+1+p$ for some integer $p\geq 0$, then it is arithmetically Cohen-Macaulay and $S_C$ satisfies property $N_{2,p}$. Recently, in the paper \cite{ein2020singularities}, L. Ein, W. Niu and J. Park have enhanced the $(2g+1+p)$-theorem for higher secant varieties as follows: if a linearly normal smooth curve $X=C\subseteq\mathbb P^r$ of degree $d$ and genus $g$ satisfies
$$
d\geq 2g+2q+p-1
$$
with $p\geq 0$, then $S^q(C)\subseteq\mathbb P^r$ is arithmetically Cohen–Macaulay and $S_{S^q(C)}$ satisfies property $N_{q+1,p}$. The paper  also provides the information on other syzygetic properties and singularities of higher secant varieties to linearly normal smooth curves of sufficiently large degree.

\item Along the same line, letting $C\subseteq\mathbb P^r$ be a linearly normal smooth curve of high degree, one may focus on the length of $(q+1)$-strand of $S^q(C)$. For the case $q=1$, the {\it gonality conjecture} due to Green-Lazarsfeld, proved by Ein-Lazarsfeld in \cite{ein2015gonality}, asserts that $K_{p,1}(C)\neq 0$ if $1\leq p\leq e-\text{gon}(C)+1$ with $e=\text{codim}C$, where $\text{gon}(C)$ is the {\it gonality} of $C$. By referring \cite[Theorem 1.2]{sidman2009syzygies}, it is easy to see that $K_{p,q}(S^q(C))\neq 0$ if $1\leq p\leq e-\gamma_q(C)+q$ with $e=\text{codim}S^q(C)$, where $\gamma_q(C)$ is the $q$-th entry in the {\it gonality sequence} of $C$ with $\gamma_1(C)=\text{gon}(C)$. So, one might expect the following: 
\begin{center}
$K_{p,q}(S^q(C))=0$ for every $p>e-\gamma_q(C)+ q$.
\end{center} 
Settling this kind of estimate, we shall discover another matryoshka phenomenon, namely the {\it generalized gonality conjecture}.

\item It would be interesting to classify smooth projective varieties whose $q$-secant varieties are del Pezzo for some integer $q\geq 2$.
\end{enumerate}

\section*{Appendix: Partial elimination ideals}
\appendix
\newcommand\invisiblesection[1]{
  \refstepcounter{section}
  \sectionmark{#1}}
\invisiblesection{}

Let $x_0,\ldots,x_r$ be the homogeneous coordinates of $\mathbb P^r$. When $z=(1:0:\cdots:0)\in\mathbb P^r$, $S_z$ denotes the polynomial ring generated by variables $x_1,\ldots,x_r$, and $V_z$ the vector space $(S_z)_1$ of linear forms vanishing at $z$.

\begin{Def}[Partial Elimination Ideals]
Let $I\subseteq S$ be a homogeneous ideal, and $z$ any point in $\mathbb P^r$. Taking a homogeneous coordinate system $(x_0:\cdots:x_r)$ so that $z=(1:0:\cdots:0)$, consider
$$
\widetilde{K}_i(z,I):=\{f\in I:\deg_{x_0}f\leq i\}.
$$
For $f\in \widetilde{K}_i(z,I)$, write $f=f'x_0^i+\text{(lower $x_0$-degree terms)}$. Define $K_i(z,I)$ to be the image of the map $f\mapsto f'$, or equivalently, $f\mapsto\partial^if/\partial x_0^i$. This is a homogeneous ideal of $S_z$, and called the {\it $i$-th partial elimination ideal} of $I$.
\end{Def}
We have a short exact sequence
\begin{equation}\label{associated}
\xymatrix{
0 \ar[r] & \widetilde{K}_{i-1}(z,I) \ar[r] & \widetilde{K}_i(z,I) \ar[r] & K_i(z,I)(-i) \ar[r] & 0
}
\end{equation}
fitting into the following commutative diagram:
\begin{equation}\label{x0commutes}
\xymatrix{
0 \ar[r] & \widetilde{K}_{i-1}(z,I) \ar[r] \ar[d]_{\times x_0} & \widetilde{K}_i(z,I) \ar[r] \ar[d]_{\times x_0} & K_i(z,I)(-i) \ar[r] \ar_{\subseteq}[d] & 0 \\
0 \ar[r] & \widetilde{K}_i(z,I)(1) \ar[r] & \widetilde{K}_{i+1}(z,I)(1) \ar[r] & K_{i+1}(z,I)(-i) \ar[r] & 0
}
\end{equation}

The chain
$$
K_0(z,I)\subseteq K_1(z,I)\subseteq\cdots\subseteq K_i(z,I)\subseteq\cdots
$$
of ideals stabilizes at $i=s$, that is, $K_{s-1}(z,I)\subsetneq K_s(z,I)=K_{s+1}(z,I)=\cdots$. Here, the {\it stabilization number $s(z,I)$} is defined by $s(z,I)=s$, and $K_\infty(z,I)$ is defined to be $K_s(z,I)$.

The following is the geometric interpretation of partial elimination ideals:

\begin{Prop}\label{geodescription}
Let $I$ be a homogeneous ideal in $S$, and $z$ a point of $\mathbb P^r$. For the projection $\pi_z:V(I)\setminus\{z\}\to\mathbb P^{r-1}$,
$$
V(K_i(z,I))=\{p\in\mathbb P^{r-1}:\textup{length}(\pi_z^{-1}(p))>i\}\cup V(K_\infty(z,I))
$$
set-theoretically.
\end{Prop}

\begin{proof}
Let $K_i=K_i(z,I)$ and $K_\infty=K_\infty(z,I)$. Write
$$
U_i=\{p\in\mathbb P^{r-1}:\textup{length}(\pi_z^{-1}(p))\leq i\}\setminus V(K_\infty).
$$
It is easy to show that $\mathbb P^{r-1}\setminus V(K_i)\subseteq U_i$. For the reverse containment, suppose that $p\in U_i$. Let $l=\text{length}(\pi_z^{-1}(p))=\text{length}((V(I)\setminus\{z\})\cap\langle z,p\rangle)$. Since $p\not\in V(K_\infty)$, there is a polynomial $f\in I$ with $\deg_{x_0}f=s$ such that $\deg_{x_0}f(x_0,p)=s$, where $s=s(z,I)$. And by Serre vanishing, there is a polynomial $g\in I$, say $\deg_{x_0}g=d$, such that $\deg_{x_0}g(x_0,p)=l$. Write $f=\sum f_ix_0^i$ and $g=\sum g_ix_0^i$ with $f_i,g_i$ not having $x_0$. Consider a Sylvester type square matrix $M=
\begin{pmatrix}
M' \\
M''
\end{pmatrix}
$, where
$$
M'=
\begin{pmatrix}
f_s & \cdots & f_l & f_{l-1} & 0  & \cdots  \\
0 & f_s & \cdots & f_l & f_{l-1} & \ddots \\
\vdots & \ddots & \ddots & \ddots & \ddots & \ddots \\
 & \cdots & 0 & f_s & \cdots & f_l
\end{pmatrix}
$$
has size $(d-l)\times(d+s-2l)$ and
$$
M''=
\begin{pmatrix}
g_d & \cdots & g_l & g_{l-1} & 0 & \cdots \\
0 & g_d & \cdots & g_l & g_{l-1} & \ddots \\
\vdots & \ddots & \ddots & \ddots & \ddots & \ddots \\
 & \cdots & 0 & g_d & \cdots &  g_l
\end{pmatrix}
$$
has size $(s-l)\times(d+s-2l)$. Then, using the classical adjoint matrix of $M$, one may find a row vector $v=(a_{d-l-1},\ldots,a_0,b_{s-l-1},\ldots,b_0)$ of polynomials such that $vM=(0,\ldots,0,\det M)$. Put $a=\sum a_ix_0^i$ and $b=\sum b_ix_0^i$ so that the polynomial $h:=af+bg\in I$ has form $(\det M)x_0^l+\text{(lower $x_0$-degree terms)}$. However, $(\det M)(p)=f_s(p)^{d-l}g_l(p)^{s-l}\neq 0$. Hence, $p\not\in V(K_l)\supseteq V(K_i)$. We are done.
\end{proof}

\begin{Cor}\label{sprime}
For a generically finite projection $\pi_z:X\dashrightarrow X_z$, we have
$$
\deg \pi_z=\min\{i:\dim V(K_i(z,I_X))<\dim X_z\}.
$$
\end{Cor}

In order to use partial elimination ideals, the following mapping cone sequence is fundamental:
\begin{equation}\label{mappingcone}
\resizebox{0.93\textwidth}{!}{$
\cdots\to K_{i,j-1}(I,V_z)\to K_{i,j}(I,V_z)\to K_{i,j}(I,V)\to K_{i-1,j}(I,V_z)\to K_{i-1,j+1}(I,V_z)\to\cdots$}
\end{equation}
In this long exact sequence, the maps $K_{i,j-1}(I,V_z)\to K_{i,j}(I,V_z)$ come from the multiplication by $x_0$.

\subsection{Hilbert polynomials via the inner projection, and applications}

Let $M$ be a finitely generated graded $S$-module of dimension $n+1$. Consider the Hilbert polynomial
$$
P_M(m)=\sum_{i=0}^n\chi_i(M)\binom{m+i-1}{i}.
$$
When $M=S_X$, we write $P_X(m)$ and $\chi_i(X)$ instead. If $X$ has dimension $n$, then $\chi_n(X)=\deg X$ and $\chi_{n-1}(X)=1-\pi(X)$.

Let $z$ be a point in $X$ such that the projection $\pi_z:X\dashrightarrow X_z$ has degree $s'<\infty$, and $s$ the stabilization number $s(z,I_X)$. By the choice of $s'$, one may find that
$$
K_0(z,I_X)=\cdots=K_{s'-1}(z,I_X).
$$
Using the short exact sequences \eqref{associated}, we have
\begin{align*}
H_{S_X}(m) & =\sum_{i=0}^\infty H_{S_z/K_i}(m-i) \\
& =\sum_{i=0}^{s'-1}H_{S_z/K_i}(m-i)+\sum_{i=s'}^{s-1}H_{S_z/K_i}(m-i)+\sum_{i=s}^\infty H_{S_z/K_\infty}(m-i) \\
& =\sum_{i=0}^{s'-1}H_{S_z/K_0}(m-i)+\sum_{i=s'}^{s-1}H_{S_z/K_i}(m-i)+H_{S/(K_\infty)}(m-s),
\end{align*}
hence
$$
P_X(m)=\sum_{i=0}^{s'-1}P_{X_z}(m-i)+\sum_{i=s'}^{s-1}P_{Z_i}(m-i)+P_{\mathbb TC_zX}(m-s),
$$
where $Z_i\subseteq\mathbb P^{r-1}$ is the projective scheme defined by $K_i(z,I_X)$. Notice that the difference $P_X(m)-P_X(m-1)$ is the Hilbert polynomial of a general hyperplane section of $X$. So, we have
$$
\begin{aligned}
(\deg X)m+1-\pi(X)& =\sum_{i=0}^{s'-1}((\deg X_z)(m-i)+1-\pi(X_z))+\sum_{i=s'}^{s-1}\deg Z_i \\
&\phantom{=}+(\deg\mathbb TC_zX)(m-s)+1-\pi(\mathbb TC_zX).
\end{aligned}
$$
Thus, the following is obtained:

\begin{Thm}\label{secgenus}
Let $z\in X$ be a point such that the projection $\pi_z:X\dashrightarrow X_z$ is generically finite. Write $s=s(z,I)$ and $s'=\deg\pi_z$. Let $Z_i\subseteq\mathbb P^{r-1}$ be the projective scheme defined by $K_i(z,I_X)$. Then, we have
\begin{enumerate}[\normalfont (1)]
\item
$$
\deg X=s'\deg X_z+\deg \mathbb TC_zX;\text{ and}
$$

\item
$$
\pi(X)=s'\pi(X_z)-s'+\binom{s'}{2}\deg X_z-\sum_{i=s'}^{s-1}\deg Z_i+\pi(\mathbb TC_zX)+s\deg \mathbb TC_zX.
$$
\end{enumerate}
\end{Thm}

To compute the Betti table, an elementary step is to calculate alternating sums of Betti numbers, that is,
$$
B_j(M,V):=\sum_i(-1)^i\beta_{i,j-i}(M)=\sum_i(-1)^i\dim K_{i,j-i}(M,V).
$$
If no confusion arises, then we simply write $B_j(M)$ instead of $B_j(M,V)$. The following is used for the syzygetic characterization of del Pezzo higher secant varieties:

\begin{Prop}\label{alternating}
Let $I\subseteq S$ be a homogeneous ideal, and $z\in V(I)$ a point with $s=s(z,I)$. Write $K_i=K_i(z,I)$. Then, for every $j$, we have
$$
B_j(I,V)=B_j(K_0,V_z)+\sum_{i=1}^s(B_{j-i}(K_i,V_z)-B_{j-i}(K_{i-1},V_z)).
$$
Especially, when $s(z,I)=1$, we have $B_j(I,V)=B_j(K_0,V_z)-B_{j-1}(K_0,V_z)+B_{j-1}(K_\infty,V_z)$ for all $j$.
\end{Prop}

\begin{proof}
By the mapping cone sequence above, it is easy to get
$B_j(I,V)=B_j(I,V_z)-B_{j-1}(I,V_z)$. Also, by letting $\widetilde{K}_i=\widetilde{K}_i(z,I)$, a sufficiently large $N$($\geq j$) induces $B_{j'}(\widetilde{K}_N)=B_{j'}(I,V_z)$ with $j'=j$ or $j-1$. And by definition,
$$
B_{j'}(\widetilde{K}_N)=B_{j'}(\widetilde{K}_{N-1})+B_{j-i}(K_N)=\cdots=\sum_{k=0}^NB_{j'-k}(K_k).
$$
So, $B_j(I,V)=\sum_{i=0}^\infty(B_{j-i}(K_i)-B_{j-i-1}(K_i))$.
\end{proof}

\subsection{A basic inequality of Betti numbers}

\begin{Prop}\label{basicA}
Let $I\subseteq S=k[x_0,\ldots,x_r]$ be a homogeneous ideal, and $z$ any point in $\mathbb P^r$. Suppose that $I_q=0$, $I_{q+1}=\widetilde{K}_1(z,I)_{q+1}$ and $I_{q+2}=\widetilde{K}_2(z,I)_{q+2}$. Then, we have
$$
\beta_{p,q}(S/I)\leq\beta_{p,q}(S_z/K_0(z,I))+\beta_{p-1,q}(S_z/K_0(z,I))+\beta_{p,q-1}(S_z/K_1(z,I)).
$$
And if $\beta_{p-1,q+1}(S_z/K_0(z,I))=0$ and $\beta_{p-1,q}(S_z/K_1(z,I))=0$, then the equality holds.
\end{Prop}

\begin{proof}
Write $\widetilde{K}_i=\widetilde{K}_i(z,I)$ and $K_i=K_i(z,I)$. Note that by assumption, we have
\begin{enumerate}[\normalfont (1)]
\item $K_{p-1,q}(I,V_z)=0$;
\item $K_{i,q+1}(I,V_z)=K_{i,q+1}(\widetilde{K}_1)$;
\item $K_{p-2,q+2}(I,V_z)=K_{p-2,q+2}(\widetilde{K}_2)$;
\item $K_{i,q-1}(K_1)=0$; and
\item $K_{p-1,q-1}(K_2/K_1)=0$.
\end{enumerate}
We begin with the following diagram induced by \eqref{mappingcone}:
$$
\resizebox{\textwidth}{!}{
\xymatrix{
K_{p-1,q}(I,V_z) \ar[r] \ar@{=}[d]& K_{p-1,q+1}(I,V_z) \ar[r] \ar@{=}[d]& K_{p-1,q+1}(I,V) \ar[r] & K_{p-2,q+1}(I,V_z) \ar[r]^{\times x_0} \ar@{=}[d] & K_{p-2,q+2}(I,V_z) \ar@{=}[d] \\
0 & K_{p-1,q+1}(\widetilde{K}_1) &  & K_{p-2,q+1}(\widetilde{K}_1) \ar[r]_{\mu}^{\times x_0} & K_{p-2,q+2}(\widetilde{K}_2).
}
}
$$
We will count the dimension of $\text{ker}\mu$. First of all, we have
$$
\dim\ker\mu=\beta_{p,q}(S/I)-\beta_{p-1,q+1}(\widetilde{K}_1).
$$
And by the exact sequence
$$
\resizebox{\textwidth}{!}{
\xymatrix{
K_{p,q-1}(K_1) \ar[r] & K_{p-1,q+1}(\widetilde{K}_0) \ar[r] & K_{p-1,q+1}(\widetilde{K}_1) \ar[r] & K_{p-1,q}(K_1) \ar[r] & K_{p-2,q+2}(\widetilde{K}_0),
}
}
$$
together with $K_{p,q-1}(K_1)=0$, we get
$$
\dim\ker\mu\geq\beta_{p,q}(S/I)-\beta_{p,q}(S_z/K_0(z,I))-\beta_{p-1,q}(K_1).
$$
Now, let us look at the commutative diagram
$$
\xymatrix{
K_{p-1,q-1}(K_1) \ar[r] & K_{p-2,q+1}(\widetilde{K}_0) \ar[r] & K_{p-2,q+1}(\widetilde{K}_1) \ar[r]^f \ar[d]_{\mu} & K_{p-2,q}(K_1) \ar[d]^{\nu}_{\times x_0} \\
 & K_{p-2,q+2}(\widetilde{K}_1) \ar[r] & K_{p-2,q+2}(\widetilde{K}_2) \ar[r]_g & K_{p-2,q}(K_2).
}
$$
Since $K_{p-1,q-1}(K_2/K_1)=0$, the map $\nu$ is injective. And as $K_{p-1,q-1}(K_1)=0$,
$$
\dim\ker\mu\leq\dim\ker(g\circ\mu)=\dim\ker(\nu\circ f)=\dim\ker f=\beta_{p-1,q}(X_z).
$$
Combine the inequalities above. Note that if
$$
\beta_{p-1,q+1}(S_z/K_0(z,I))=0\quad\text{and}\quad\beta_{p-1,q}(S_z/K_1(z,I))=0,
$$
hence $K_{p-2,q+2}(\widetilde{K}_1)=0$, then the equality holds.
\end{proof}

\subsection{A vanishing theorem}

\begin{Lem}
Let $I\subseteq S$ be a homogeneous ideal, and $z$ any point in $\mathbb P^r$. Write $\widetilde{K}_i=\widetilde{K}_i(z,I)$ and $K_i=K_i(z,I)$. Suppose that $s(z,I)=1$. Then, the map $\times x_0:K_{i,j}(I,V_z)\to K_{i,j+1}(I,V_z)$ is
\begin{enumerate}[\normalfont (1)]
\item injective if $K_{i,j}(K_0,V_z)\to K_{i,j}(K_1,V_z)$ is one-to-one; and

\item surjective if
\begin{enumerate}[\normalfont (a)]
\item $K_{i,j}(K_0,V_z)\to K_{i,j}(K_1,V_z)$ is onto;

\item $K_{i+1,j-1}(K_1,V_z)\to K_{i,j+1}(\widetilde{K}_0,V_z)$ is onto; and

\item $K_{i-1,j+1}(K_0,V_z)\to K_{i-1,j+1}(K_1,V_z)$ is one-to-one.
\end{enumerate}
\end{enumerate}
\end{Lem}

\begin{proof}
For a sufficiently large $N$($\geq j$), due to the commutative diagram
$$
\xymatrix{
K_{i,j}(\widetilde{K}_N) \ar[r]^--{\times x_0} \ar[d]_\cong & K_{i,j+1}(\widetilde{K}_{N+1}) \ar[d]^\cong \\
K_{i,j}(I,V_z) \ar[r]_--{\times x_0} & K_{i,j+1}(I,V_z),
}
$$
it is enough to consider the map
$$
\times x_0:K_{i,j}(\widetilde{K}_N)\to K_{i,j+1}(\widetilde{K}_{N+1}).
$$
We will use the commutating diagram \eqref{x0commutes}. Note that the map $K_k(-k)\to K_{k+1}(-k)$ is an isomorphism whenever $k\geq 1$ because $s(z,I)=1$. Taking the Koszul cohomology, we have
$$
\resizebox{\hsize}{!}{
\xymatrix{
K_{i+1,j-k-1}(K_k) \ar[r] \ar[d] & K_{i,j}(\widetilde{K}_{k-1}) \ar[r] \ar[d]_{\times x_0} & K_{i,j}(\widetilde{K}_k) \ar[r] \ar[d]_{\times x_0} & K_{i,j-k}(K_k) \ar[r] \ar[d] & K_{i-1,j+1}(\widetilde{K}_{k-1}) \ar[d]_{\times x_0} \\
K_{i+1,j-k-1}(K_{k+1}) \ar[r] & K_{i,j+1}(\widetilde{K}_k) \ar[r] & K_{i,j+1}(\widetilde{K}_{k+1}) \ar[r] & K_{i,j-k}(K_{k+1}) \ar[r] & K_{i-1,j+2}(\widetilde{K}_k).
}
}
$$
Setting $k=0$, we have
$$
\xymatrix{
& & K_{i,j}(\widetilde{K}_0) \ar[r]^\cong \ar[d]_{\times x_0} & K_{i,j}(K_0) \ar[d] \\
K_{i+1,j-1}(K_1) \ar[r] &K_{i,j+1}(\widetilde{K}_0) \ar[r] & K_{i,j+1}(\widetilde{K}_1) \ar[r] & K_{i,j}(K_1).
}
$$

We claim that
\begin{enumerate}[\normalfont (1)]
\item[(1)] $\times x_0:K_{i,j}(\widetilde{K}_k)\to K_{i,j+1}(\widetilde{K}_{k+1})$ is injective for all $k\geq 0$ provided that $K_{i,j}(K_0)\to K_{i,j}(K_1)$ is a monomorphism.

\item[(2)] $\times x_0:K_{i,j}(\widetilde{K}_k)\to K_{i,j+1}(\widetilde{K}_{k+1})$ is surjective for all $k\geq 0$ provided that $K_{i,j}(K_0)\to K_{i,j}(K_1)$ and $K_{i+1,j-1}(K_1)\to K_{i,j+1}(\widetilde{K}_0)$ are epimorphisms, and $K_{i-1,j+1}(K_0)\to K_{i-1,j+1}(K_1)$ is injective.
\end{enumerate}

We find that $\times x_0:K_{i,j}(\widetilde{K}_0)\to K_{i,j+1}(\widetilde{K}_1)$ is an injection when so is $K_{i,j}(K_0)\to K_{i,j}(K_1)$,  Applying the $4$-lemma, we have (1).

For the proof of (2), we will use induction on $k\geq 0$. The case $k=0$ is obvious by the surjectivity of $K_{i,j}(K_0)\to K_{i,j}(K_1)$ and $K_{i+1,j-1}(K_1)\to K_{i,j+1}(K_0)$. Set $k\geq 1$.  Since $K_{i-1,j+1}(K_0)\to K_{i-1,j+1}(K_1)$ is injective, so is $\times x_0:K_{i-1,j+1}(\widetilde{K}_k)\to K_{i-1,j+2}(\widetilde{K}_{k+1})$ by (1). Then, by the $4$-lemma on surjectivity, we are done.
\end{proof}

\begin{Prop}\label{vanishingA}
Suppose that $s(z,I)=1$. Then, $K_{i,j}(I,V)=0$ if
\begin{enumerate}[\normalfont (1)]
\item $K_{i,j-1}(K_\infty,V_z)=0$; and

\item $K_{i-1,j}(K_0,V_z)=K_{i,j}(K_0,V_z)=0$.
\end{enumerate}
\end{Prop}

\begin{proof}
Consider the exact sequence \eqref{mappingcone}, that is,
$$
\resizebox{\textwidth}{!}{
\xymatrix{
K_{i,j-1}(I,V_z) \ar[r]^--{\times x_0} & K_{i,j}(I,V_z) \ar[r] & K_{i,j}(I,V) \ar[r] & K_{i-1,j}(I,V_z) \ar[r]^--{\times x_0} & K_{i-1,j+1}(I,V_z).
}
}
$$
By applying the lemma above, the first $\times x_0$ is surjective, and the second $\times x_0$ is injective.
\end{proof}

\nocite{*}
\bibliographystyle{amsplain}
\bibliography{matryoshka_202110}

\end{document}